\documentclass[aop,seceqn,MSNbibl,citesort,dvips]{arximspdf}
\usepackage{graphicx}

\doi{10.1214/10-AOP572}
\volume{39}
\issue{3}
\pubyear{2011}
\firstpage{1161}
\lastpage{1203}

\makeatletter

\newcommand{\ra}{\rightarrow}

\newtheorem{theo}{Theorem}[section]
\newtheorem{pro}[theo]{Proposition}
\newtheorem{lem}[theo]{Lemma}
\newtheorem{cor}[theo]{Corollary}
\newtheorem{claim}[theo]{Claim}

\newproclaim{defin}[theo]{Definition}
\newproclaim{prob}[theo]{Problem}
\newproclaim{rem}[theo]{Remark}

\makeatother

\begin{document}
\begin{frontmatter}

\title{Merging for inhomogeneous finite Markov chains, part II: Nash
and log-Sobolev
inequalities}
\runtitle{Merging for inhomogeneous Markov chains}

\begin{aug}
\author[A]{\fnms{L.} \snm{Saloff-Coste}\thanksref{t1}} and
\author[B]{\fnms{J.} \snm{Z\'u\~{n}iga}\corref{}\thanksref{t2}
\ead[label=e1]{jzuniga@math.stanford.edu}}
\runauthor{L. Saloff-Coste and J. Z\'u\~{n}iga}
\affiliation{Cornell University and Stanford University}
\address[A]{Department of Mathematics\\
Cornell University\\
Malott Hall\\
Ithaca, New York, 14853\\
USA} %adresu isvedimo komanda gale!
\address[B]{Department of Mathematics\\
Stanford University\\
Stanford, California, 94305\\
USA\\
\printead{e1}}
\end{aug}

\thankstext{t1}{Supported in part by NSF Grant DMS-06-03886.}

\thankstext{t2}{Supported in part by NSF Grants DMS-06-03886,
DMS-03-06194 and by a NSF Postdoctoral Fellowship.}

% HISTORY:
\received{\smonth{3} \syear{2009}}
\revised{\smonth{4} \syear{2010}}

% ABSTRACT
%
\begin{abstract}
We study time-inhomogeneous Markov chains with finite state spaces
using Nash and logarithmic-Sobolev inequalities, and the notion of
$c$-stability. We develop the basic theory of such functional
inequalities in the time-inhomogeneous context and provide illustrating
examples.
\end{abstract}

% KEYWORDS
%
\setattribute{keyword}{AMS}{MSC2010 subject classification.}

\end{frontmatter}

%s1 ###
\section{Introduction}

%s1.1 ###
\subsection{Background}
This article is part of a series of works where we study quantitative
merging properties of time inhomogeneous finite Markov chains.
Time inhomogeneity leads to a great variety of behaviors.
Moreover, even in rather simple situations, we are at a loss to study
how a time inhomogeneous Markov chain might behave. Here, we focus on
a natural but restricted type of problem. Consider a sequence
of aperiodic irreducible Markov kernels $(K_i)_1^{\infty}$ on a finite
set $V$.
Let $\pi_i$ be the invariant measure of $K_i$. Assume that, in a sense to
be made precise, all $K_i$ and all $\pi_i$ are similar and the behavior
of the time homogeneous chains driven by each $K_i$ separately is understood.
Can we then describe the behavior of the time inhomogeneous chain
driven by the sequence $(K_i)_1^{\infty}$?

To give a concrete example, on
$V_N=\{0,\ldots,N\}$, consider a sequence of aperiodic irreducible
birth and death chain kernels $K_i$, $i=1,2,\ldots,$ with
\[
1/4\le K_i(x,y)\le3/4 \qquad\mbox{if } |x-y|\le1
\]
and with
reversible measure $\pi_i$ satisfying
$1/4\le(N+1)\pi_i(x)\le4$, for all $x\in V_N$.
What can we say about the behavior of
the corresponding time inhomogeneous Markov chain?

Remarkably enough, there is very little known about
this question. What can we expect to be true?
What can we try to prove? Let $K_{0,n}(x,\cdot)$ denote
the distribution, after $n$ steps,
of the time inhomogeneous chain described above
started at $x$. It is not hard to see that such a chain
satisfies a Doeblin type condition that implies
\[
{\lim_{n\ra\infty}}\|K_{0,n}(x,\cdot)-K_{0,n}(y,\cdot)
\|_{\mathrm{TV}}=0.
\]
In the absence of a true target distribution and following
\cite{ADF}, we call this property \textit{merging}.
Of course, this does not qualify as a quantitative result.
Extrapolating from the behavior of each kernel $K_i$ taken
individually, we may hope to show that, if
$ \lim_{N\ra\infty} t_N/N^2=\infty$ then
\[
{\lim_{N\ra\infty}}\|K_{0,t_N}(x,\cdot)-K_{0,t_N}(y,\cdot)
\|_{\mathrm{TV}}=0.
\]

The aim of this paper and the companion paper \cite{SZ3} is to present
techniques that apply to this
type of problem. The simple minded problem outlined above is
actually quite challenging and we will not be able to resolve it
here without some additional hypotheses. However, we show how to
adapt techniques such as singular values, Nash and log-Sobolev inequalities
to time inhomogeneous chains and provide
a variety of examples where these tools apply.
In \cite{SZ3}, we discussed singular value techniques.
Here, we focus on Nash and log-Sobolev inequalities.
The examples treated here (as well as those treated in \cite{SZ3,SZ-wave})
are quite particular despite the fact that one may believe that
the techniques we use are widely applicable. Whether or not
such a belief
is warranted is a very interesting and, so far, unanswered question.
This is deeply related to the notion of $c$-stability that is introduced
here and in \cite{SZ3}. The examples we present here and in
\cite{SZ,SZ3,SZ-wave} are about the only existing
evidence of successful quantitative analysis of time inhomogeneous
Markov chains.

A more detailed introduction to these questions is in \cite{SZ3}.
The references
\cite{Ga,SZ} discuss singular value techniques in the case of time
inhomogeneous chains that admit an invariant distribution [all kernels
$K_i$ in the sequence $(K_i)_1^\infty$ share a common invariant distribution].
Time inhomogeneous random walks on finite groups provide a large collection
of such
examples (see also \cite{MPS} for a particularly interesting example:
semirandom transpositions). The papers \cite{DLM,DMR} are also concerned
with quantitative results for time
inhomogeneous Markov chains. In particular,
the techniques developed in \cite{DLM} are closely related
to ours and we will use some of their results concerning
the modified logarithmic Sobolev inequality.
References on the basic theory of
time inhomogeneous Markov chains are \cite{Io,Pa,Sen,Sen2,Son}. For a different
perspective, see also \cite{CH}.

A short review of the relevant aspects of the time inhomogeneous Markov chain
literature, including the use of ``ergodic coefficients'' can be found
in \cite{SZ-berlin}.
The vast literature on the famous simulated annealing algorithm is not
very relevant
for our purpose but we refer to \cite{DM} for a recent discussion. The
paper \cite{DG}
concerned with filtering and genetic algorithms describes problems that
are related in
spirit to the present work.

%s1.2 ###
\subsection{Basic notation}\label{subsec-BN}
Let $V$ be a finite set equipped with a sequence of kernels
$(K_n)_1^\infty$ such that, for each $n$,
$K_n(x,y)\ge0$ and $\sum_y K_n(x,y)=1$. An associated Markov chain
is a $V$-valued random process $X=(X_n)_0^\infty$ such that, for all $n$,
\begin{eqnarray*}
P(X_n=y|X_{n-1}=x,\ldots,X_0=x_0)&=&P(X_n=y|X_{n-1}=x)\\
&=&
K_n(x,y).
\end{eqnarray*}
The distribution $\mu_n$ of $X_n$ is determined by
the initial distribution
$\mu_0$ and given by
\[
\mu_n(y)=\sum_{x\in V}\mu_0(x)K_{0,n}(x,y),
\]
where $K_{n,m}(x,y)$ is defined inductively for each $n$ and each $m\ge n$
by
\[
K_{n,m}(x,y)=\sum_{z\in V} K_{n,m-1}(x,z)K_m(z,y)
\]
with $K_{n,n}=I$ (the identity).
If we interpret the $K_n$'s as matrices, then this definition means
that $K_{n,m}=K_{n+1}\cdots K_m$. This paper is mostly
concerned with the behavior
of the measures $K_{0,n}(x,\cdot)$ as $n$ tends to infinity.
In the case of time homogeneous chains where all $K_i=Q$ are equal, we write
$K_{0,n}=Q^n$.

Our main interest is in ergodic like properties of time
inhomogeneous Markov chains. In general, one does not expect
$\mu_n=\mu_0K_{0,n}$ to converge toward a limiting distribution.
Instead, the natural notion is that of merging of measures as discussed in
\cite{ADF}.
\begin{defin} Fix a sequence of Markov kernels as above.
We say the sequence
is merging if for any $x,y,z\in V$,
%
%e1.1 ###
%
\begin{equation}\label{TVmerg}
\lim_{n\ra\infty}K_{0,n}(x,z)-K_{0,n}(y,z)=0.
\end{equation}
\end{defin}
\begin{rem} If the sequence $(K_i)_1^\infty$ is merging then, for any
two starting distributions $\mu_0,\nu_0$, the measures $\mu_n=\mu_0K_{0,n}$
and $\nu_n=\nu_0K_{0,n}$ are merging, that is, $\mu_n-\nu_n\ra0$.
Since we assume the set $V$ is finite, merging is equivalent to
${\lim_{n\ra\infty}}\|K_{0,n}(x,\cdot)-K_{0,n}(y,\cdot)\|_{\mathrm{TV}}=0$.
Hence, we also refer to this property as ``total variation merging.''
\end{rem}

Total variation merging is also referred to as weak ergodicity in the
literature and there exists a body of work concerned with understanding
when weak ergodicity holds. See, for example, \cite{Io,NS,Pa,Rho,Sen}.
A main tool used to show weak ergodicity is that of contraction
coefficients. Furthermore, in \cite{FJ}, Birkhoff's contraction
coefficient is used to study \textit{ratio ergodicity} which is equivalent
to what we will later call relative-sup merging. However, it should be
noted that even for time homogeneous chains Birkhoff coefficients and
related methods fail to provide useful quantitative bounds in most cases.

Our goal is to develop quantitative results
in the context of time inhomogeneous chains in the spirit of the work of
Aldous, Diaconis and others. In these works, precise estimates of the
mixing time of ergodic chains are obtained. Typically,
a family of Markov chains indexed by a parameter, say $N$, is studied.
Loosely speaking, as the parameter $N$ increases, the complexity and
size of the chain increases and one seeks bounds that depend on $N$
in an explicit quantitative way. See, for example,
\cite{Al,BT,Dia,DSh,DS-C,DS-N,DS-L,DS-M,Fill,MT,MP,StF}.
Efforts in this direction for time inhomogeneous chains
are in \cite{DLM,DMR,FJ,Ga,Goel,MPS,SZ,SZ3}. Still, there are only a
very small number of
results and examples concerning the quantitative study of merging as
defined above
for time inhomogeneous Markov chains so that it is not very clear what
kind of results should be expected and what kind of hypotheses are
reasonable. We
refer the reader to \cite{SZ3} for a more detailed discussion.

The following definition is useful to capture the spirit of our study.
It indicates that the simplest case we would like to
think about is the case when the sequence $K_i$ is obtained by
deterministic but arbitrary choices between a finite number of kernels
$\mathcal Q=\{Q_1,\ldots,Q_k\}$.
\begin{defin}
We say that a set
$\mathcal{Q}$ of Markov kernels on
$V$ is merging in total variation if for any sequence $(K_i)_0^{\infty}$
with $K_i\in\mathcal{Q}$ for all $i$,
we have
\[
\forall x,y,z\in V\qquad
{\lim_{n\ra\infty}}\|K_{0,n}(x,\cdot)-K_{0,n}(y,\cdot)\|_{\mathrm{TV}}=0.
\]
\end{defin}

In the study of ergodicity of finite Markov chains, the convergence
toward the
target distribution is measured using various notions of distance between
probability measures. These include the total variation distance
\[
\|\mu-\nu\|_{\mathrm{TV}}=\sup_{A\subset V}\{\mu(A)-\nu(A)\},
\]
the chi-square distance (w.r.t. $\nu$. Note the asymmetry between $\mu$
and $\nu$.)
\[
\biggl(\sum_y\biggl|\frac{\mu(y)}{\nu(y)}-1\biggr|^2\nu(y)\biggr)^{1/2},
\]
and the relative sup-distance
(again, note the asymmetry)
\[
\max_y \biggl\{\biggl|\frac{\mu(y)}{\nu(y)} -1\biggr|\biggr\}.
\]
These will be used here to measure merging.

%s1.3 ###
\subsection{Merging time} In the quantitative theory of ergodic time
homogeneous Markov chains, the notion of mixing time plays a crucial role.
For time inhomogeneous chain, we propose to consider the following
definitions.
\begin{defin}
Fix $\varepsilon\in(0,1)$.
Given a sequence $(K_i)_1^\infty$ of Markov kernels on a finite
set $V$, we call max total variation merging time the quantity
\[
T_{\mathrm{TV}}(\varepsilon)=
\inf\Bigl\{n\dvtx{\max_{x,y\in V}}
\|K_{0,n}(x,\cdot)-K_{0,n}(y,\cdot)\|_{\mathrm{TV}}<\varepsilon\Bigr\}.
\]
\end{defin}
\begin{defin} Fix $\varepsilon\in(0,1)$.
We say that a set
$\mathcal{Q}$ of Markov kernels
on $V$ has max total variation $\varepsilon$-merging time at most $T$
if for any sequence $(K_i)_1^{\infty}$
with $K_i\in\mathcal{Q}$ for all $i$,
we have $T_{\mathrm{TV}}(\varepsilon)\le T$, that is,
\[
\forall t>T\qquad
\max_{x,y\in V}\{\|K_{0,t}(x,\cdot)-K_{0,t}(y,\cdot)\|_{\mathrm{TV}}
\}
\le\varepsilon.
\]
\end{defin}

Of course, merging can be measured in ways other than total variation. Also
merging is a bit less flexible than mixing in this respect since there
is no reference measure. One very natural and much stronger notion
than total variation is relative sup-distance. For time inhomogeneous chains,
total variation merging does not necessarily
imply relative-sup merging as defined below. See \cite{SZ3}.
\begin{defin}
We say a sequence $(K_i)_1^\infty$ of Markov kernels on a finite
set $V$ is merging in relative-sup if for all $x,y,z\in V$
\[
\lim_{n\ra\infty}\frac{K_{0,n}(x,z)}{K_{0,n}(y,z)}=1
\]
with the convention that $0/0=1$ and $a/0=\infty$ for $a>0$.
Fix $\varepsilon\in(0,1)$, we call relative-sup merging time the quantity
\[
T_{\infty}(\varepsilon)=
\inf\biggl\{ n\dvtx\max_{x,y,z\in V}
\biggl\{\biggl|\frac{K_{0,n}(x,z)}{K_{0,n}(y,z)}-1\biggr|\biggr\}
<\varepsilon\biggr\}.
\]
\end{defin}
\begin{defin}
We say a set $\mathcal{Q}$ of Markov kernels on $V$ is merging in
relative-sup if any
sequence $(K_i)_1^{\infty}$ with $K_i\in\mathcal Q$ for all $i$ is
merging in relative-sup.

Fix $\varepsilon\in(0,1)$.
We say that $\mathcal{Q}$ has relative-sup $\varepsilon$-merging time at
most $T$
if for any sequence $(K_i)_1^{\infty}$
with $K_i\in\mathcal{Q}$ for all $i$,
we have $T_{\infty}(\varepsilon)\le T$, that is,
\[
\forall t>T\qquad
\max_{x,y,z\in V}\biggl\{\biggl|\frac{K_{0,t}(x,z)}{K_{0,t}(y,z)}-1
\biggr|\biggr\}
\le\varepsilon.
\]
\end{defin}

The following problem is open. It is a quantitative version of
the problem stated at the beginning of the introduction.
\begin{prob} \label{Pb-bdu}
Let $V_N=\{0,\ldots,N\}$ and $c\in[1,\infty)$.
Let $\mathcal Q_N$ be the set of all
birth and death chains $Q$ on $V_N$ with $Q(x,y)\in[1/4,3/4]$
if $|x-y|\le1$,
and reversible measure $\pi$ satisfying $1/4\le(N+1)\pi(x)\le4$,
$x\in V_N$.
\begin{enumerate}
\item Prove or disprove that there exists
a constant $A$ independent of $N$ such that
$\mathcal Q_N$ has total variation $\varepsilon$-merging time at most
$AN^2(1+\log_+ 1/\varepsilon)$.
\item
Prove or disprove that there exists
a constant $A$ independent of $N$ such that
$\mathcal Q_N$ has relative-sup $\varepsilon$-merging time at most
$AN^2(1+\log_+ 1/\varepsilon) $.
\end{enumerate}
\end{prob}
\begin{rem} This problem is open (in most cases)
even if one considers a sequence $(K_i)_1^{\infty}$
drawn from a set $\mathcal Q=\{K_1,K_2\}$ of two kernels. Observe that
the hypothesis that the invariant measures $\pi_i$ are all comparable
to the uniform plays some role.
How to harvest the global hypothesis of comparable
stationary distributions $\pi_i$ is not entirely clear. See Theorem \ref{th-bdu}
below for a partial solution.

If $\pi_1$ and $\pi_2$ are not comparable, it is possible for
$(K_1,\pi_1)$ and $(K_2,\pi_2)$ to have the same mixing time yet for
$\mathcal Q=\{K_1,K_2\}$ to have a merging time of a higher order.
Assume that $K_1$ and $K_2$ are two biased random
walks with equal drift, one drift to left, the other to the right.
Despite the
fact that each of these random walks has a relative-sup mixing time
of order $N$, the inhomogeneous chain driven by the sequence
$K_1K_2K_1K_2\cdots$ has a relative-sup merging time of order $N^2$, see
\cite{SZ3}.
%Showing that $T_{\mbox{{TV}}}(\varepsilon)\leq CN^2$ for an
%arbitrary sequence drawn from $\{K_1,K_2\}$ is an interesting open
%problem.
\end{rem}

%s1.4 ###
\subsection{Stability}
In this section, we consider a property, $c$-stability, that plays a
crucial role in the techniques
we develop to provide quantitative bounds for time inhomogeneous Markov
chains. This
property was introduced and discussed in~\cite{SZ3}. It is a straightforward
generalization of the property of sharing the same invariant measure.
Unfortunately,
it is hard to check.
\begin{defin} Fix $c\ge1$.
A sequence of Markov kernels $(K_n)_1^{\infty}$ on a finite set $V$
is $c$-stable if there exists a measure $\mu_0$ such that
%
%e1.2 ###
%
\begin{equation}\label{c-stab-seq}
\forall n\ge0, x\in V\qquad c^{-1}\leq\frac{\mu_n(x)}{\mu
_0(x)}\leq c,
\end{equation}
where $\mu_n=\mu_0K_{0,n}$. If this holds, we say that
$(K_n)_1^{\infty}$
is $c$-stable with respect to the measure $\mu_0$.
\end{defin}
\begin{defin} A set $\mathcal{Q}$ of Markov kernels
is $c$-stable with respect to a measure $\mu_0$ if any sequence
$(K_i)_1^{\infty}$ such that $K_i\in\mathcal{Q}$ for all $i$ is
$c$-stable with respect to $\mu_0$.
\end{defin}
\begin{rem} If all $K_i$ share the same invariant distribution $\pi$
then $(K_i)_1^\infty$ is $1$-stable with respect to $\pi$.
\end{rem}
\begin{rem} Suppose a set
$\mathcal{Q}$ of aperiodic irreducible Markov kernels
is $c$-stable with respect to a measure $\mu_0$. Let
$\pi$ be an invariant measure for some $Q\in\mathcal Q$.
Then we must have
\[
x\in V,\qquad \frac{1}{c}\leq\frac{\pi(x)}{\mu_0(x)}\leq c.
\]
Hence, $\mathcal{Q}$ is also $c^2$-stable with respect to $\pi$
and any two invariant measures $\pi,\pi'$ for kernels $Q,Q'\in\mathcal Q$
must satisfy
\[
x\in V,\qquad \frac{1}{c^2}\leq\frac{\pi(x)}{\pi'(x)}\leq c^2.
\]
\end{rem}

The following theorem which relates to a special case of Problem
\ref{Pb-bdu}
illustrates the role of $c$-stability.
\begin{theo}\label{th-bdu}
Let $V_N=\{0,\ldots,N\}$.
Let $\mathcal Q_N$ be the set of all
birth and death chains $Q$ on $V_N$ with
\[
Q(x,y)\in[1/4,3/4] \qquad\mbox{if } |x-y|\le1
\]
and reversible measure $\pi$ satisfying $1/4\le(N+1)\pi(x)\le4$,
$x\in V_N$. Let $(K_i)_1^\infty$ be a sequence of birth and death
Markov kernels on $V_N$ with $K_i\in\mathcal Q_N$.
Assume that $(K_i)_1^\infty$ is $c$-stable with respect to
the uniform measure on $V_N$, for some constant $c\ge1$ independent of $N$.
Then there exists a constant $A=A(c)$ (in particular, independent of $N$)
such that the relative-sup merging time for $(K_i)_1^\infty$
on $V_N$ is bounded by
\[
T_{\infty}(\varepsilon)\le AN^2(1+\log_+ 1/\varepsilon).
\]
\end{theo}

This will be proved later in a stronger form in Section \ref{section-Nash}.
In \cite{SZ3}
the weaker conclusion $T_\infty(\varepsilon)\le AN^2(\log N+\log_+
1/\varepsilon)$
was obtained using singular value techniques. Here,
we will use Nash inequalities to obtain
$T_\infty(\varepsilon)\le AN^2(1+\log_+ 1/\varepsilon)$.

It is possible that the set $\mathcal Q_N$ is $c$-stable with respect
to the
uniform measure for some $c$. Indeed, it is tempting to conjecture that
this is the case although the evidence is rather limited (see also the
discussion in \cite{SZ-berlin}).
If this is true, then Theorem \ref{th-bdu}
solves Problem \ref{Pb-bdu}. However, we do not know how to approach the
problem of proving $c$-stability for $\mathcal Q_N$.
\begin{rem}
While the assumption of $c$-stability in Theorem \ref{th-bdu}
is quite strong, Sections 4.2 and 5 of \cite{SZ3} give specific
examples of families $\mathcal Q_N$ for which it holds. Further, we
note that the question of whether or not $c$-stability holds is
extremely natural and interesting in itself.
\end{rem}

%s2 ###
\section{Singular values and Nash inequalities}

One key idea in the study of Markov chains is to associate to a Markov
kernel $K$
the operator $K\dvtx f\mapsto Kf=\sum_y K(\cdot,y)f(y)$. In the case of time
homogeneous chains, one uses the basic fact that this operator acts on
$\ell^p(\pi)$ with norm $1$ when $\pi$ is an invariant measure.

In the case of time inhomogeneous chains, it is crucial to consider $K$
as an
operator between $\ell^p$ spaces with different measures in the domain
and target
spaces. The following simple observation is key.

Given a measure $\mu$ and a Markov kernel $K$ on a finite set $V$, set
$\mu'=\mu K$.
Fix $p\in[1,\infty)$ and consider $K$ as a linear
operator
%
%e2.1 ###
%
\begin{equation}\label{Kop}
K=K_\mu\dvtx\ell^p(\mu')\ra\ell^p(\mu),\qquad Kf(x)=\sum_yK(x,y)f(y).
\end{equation}
Then
%
%e2.2 ###
%
\begin{equation}\label{cont}\qquad
\|K\|_{\ell^p(\mu')\ra\ell^p(\mu)}=
\sup\bigl\{ \|Kf\|_{\ell^p(\mu)}\dvtx f\in\ell^p(\mu'),
\|f\|_{\ell^p(\mu')}\le1\bigr\} =1.
\end{equation}
This follows from Jensen's inequality. See, for example, \cite{DLM,SZ3}.
We will use the notation $K_\mu$ whenever we need to emphasize the fact that
$K$ is viewed as an operator between $\ell^p(\mu K)$ and
$\ell^q(\mu)$ for some $1\le p,q\le\infty$. When the context is clear,
we will drop the subscript $\mu$ as was done above.

%s2.1 ###
\subsection{Using various distances}
Given a sequence of Markov kernels $(K_i)_1^\infty$,
fix a starting measure $\mu_0$ and set $\mu_n=\mu_0 K_{0,n}$.
We will assume that $\mu_n>0$ for all $n$.
Note that if $\mu_0>0$ and $K_n$ are all irreducible then $\mu_n>0$ for
all $n\geq0$.
We are interested in the behavior of
\[
d_p(K_{0,n}(x,\cdot),\mu_n)=\biggl(
\sum_y
\biggl|\frac{K_{0,n}(x,y)}{\mu_n(y)}-1\biggr|^p\mu_n(y)\biggr)^{1/p},\qquad
p\geq1.
\]
For $p\geq1$, a classical argument involving the duality between
$\ell^p$ and $\ell^q$ where $1=1/p+1/q$, yields
\[
d_p(K_{0,n}(x,\cdot),\mu_n)=\sup\biggl\{
\biggl|\sum_y [K_{0,n}(x,y)f(y) -\mu_n(y)f(y)]\biggr|\dvtx
\|f\|_{\ell^q(\mu_n)}\le1 \biggr\}
\]
and one checks that the function
\[
n\mapsto d_p(K_{0,n}(x,\cdot),\mu_n)
\]
is nonincreasing (see \cite{SZ3}). Of course,
\[
2\|K_{0,n}(x,\cdot)-\mu_n\|_{\mathrm{TV}}= d_1(K_{0,n}(x,\cdot),\mu_n)
\]
and, if $1\le p\le r\le\infty$,
\[
d_p(K_{0,n}(x,\cdot),\mu_n)\le d_{r}(K_{0,n}(x,\cdot),\mu_n).
\]
In particular,
%
%e2.3 ###
%
\begin{equation}\label{tv2}
2\|K_{0,n}(x,\cdot)-\mu_n\|_{\mathrm{TV}}\le d_2(K_{0,n}(x,\cdot
),\mu_n)
\end{equation}
and
%
%e2.4 ###
%
\begin{equation}\label{tv2xy}
\|K_{0,n}(x,\cdot)-K_{0,n}(y,\cdot)
\|_{\mathrm{TV}}\le\max_{x\in V}\{d_2(K_{0,n}(x,\cdot),\mu_n)\}.
\end{equation}
Further, if
\[
\max_{x,z}\biggl\{\biggl|\frac{K_{0,n}(x,z)}{\mu_n(z)}-1\biggr|\biggr\}\le
\varepsilon\le1/2,
\]
then
\[
\max_{x,y,z}\biggl\{\biggl|\frac{K_{0,n}(x,z)}{K_{0,n}(y,z)}-1
\biggr|\biggr\}
\le4\varepsilon.
\]
To see the last inequality, note that if
$1-\varepsilon\le a/b,c/b\le1+\varepsilon$ with $\varepsilon\in(0,1/2)$ then
\[
1-2\varepsilon\le\frac{1-\varepsilon}{1+\varepsilon}\le\frac{a}{c}\le
\frac{1+\varepsilon}{1-\varepsilon}\le1+4\varepsilon.
\]

%s2.2 ###
\subsection{Singular values}\label{sec-sv}
In \cite{SZ3}, we developed basic inequalities for $d_2(K_{0,n}(x$, $\cdot
),\mu_n)$
based on singular value decompositions.
The basic fact here is that, if $\mu$ is a probability measure on $V$,
$K$ a Markov kernel and $\mu'=\mu K$ then
\[
d_{2}(K(x,\cdot),\mu')^2=\sum_{i=1}^{|V|-1}|\psi_i(x)|^2\sigma_i^2,
\]
where $\sigma_i$, $i=0,\ldots,|V|-1$, are the singular values of
$K_\mu\dvtx
\ell^2(\mu')\ra\ell^2(\mu)$ in nonincreasing order,
that is the square root of the eigenvalues of $K_\mu K^*_\mu\dvtx
\ell^2(\mu)\ra\ell^2(\mu)$ where $K^*_\mu\dvtx\ell^2(\mu)\ra\ell^2(\mu')$
is the adjoint of
$K_\mu\dvtx\ell^2(\mu')\ra\ell^2(\mu)$.
The $\psi_i$'s form an orthonormal basis
for $\ell^2(\mu)$ and are eigenfunctions of $K_\mu K^*_\mu$,
$\psi_i$ being associated with $\sigma^2_i$.
Of course, the $\sigma_i^2$'s can also be viewed
as the eigenvalues of $K_\mu^*K_\mu\dvtx\ell^2(\mu')\ra\ell^2(\mu')$.

In any case, a crucial fact for us here is that $\sigma_1$,
the second largest singular value of $K_\mu\dvtx\ell^2(\mu')\ra\ell
^2(\mu)$,
is also the norm
of $K-\mu'=K_\mu-\mu'\dvtx\ell^2(\mu')\ra\ell^2(\mu)$, that is,
\[
\sup\bigl\{ \| (K-\mu')f\|_{\ell^2(\mu)}\dvtx f\in\ell^2(\mu'),
\|f\|_{\ell^2(\mu')}=1\bigr\}= \sigma_1.
\]

Given a sequence $(K_i)_1^\infty$ of Markov kernels on $V$ and a
positive measure
$\mu_0$, set $\mu_n=\mu_0K_{0,n}$ and let\vspace*{1pt}
$\sigma_1(K_i,\mu_{i-1})$ be the second largest singular value of
$K_i\dvtx
\ell^2(\mu_{i})\ra\ell^2(\mu_{i-1})$.
Noting that
\[
(K_{0,n}-\mu_n)= (K_1-\mu_1)(K_2-\mu_2)\cdots(K_n-\mu_n),
\]
we obtain
%
%e2.5 ###
%
\begin{equation}
\|K_{0,n}-\mu_n\|_{\ell^2(\mu_n)\ra\ell^2(\mu_0)}\le
\prod_1^n\sigma_1(K_i,\mu_{i-1}).
\end{equation}

This inequality seems very promising and this is rather misleading.
There is very little hope to compute or estimate the singular values
$\sigma_i(K_i,\mu_{i-1})$, even if we have a good grasp on the kernel $K_i$.
The reason is that $\sigma_1(K_i,\mu_{i-1})$ depends very much
on the unknown measure $\mu_{i-1}$. This is similar to the problem one faces
when studying an irreducible aperiodic time homogeneous finite Markov chain
for which one is not able to compute the stationary measure (although
this case is
rarely discussed, it is the typical case). For positive examples
and a more detailed discussion, see \cite{SZ3}.

%s2.3 ###
\subsection{Dirichlet forms} Given a reversible Markov kernel $Q$ with
reversible measure $\pi$ on a finite set $V$, the associated Dirichlet
form is
\begin{eqnarray*}
\mathcal E(f,f)&=&\mathcal E_{Q,\pi}(f,f)
=\langle(I-Q)f,f\rangle_\pi\\
&=& \frac{1}{2}\sum_{x,y}|f(x)-f(y)|^2 \pi(x)Q(x,y).
\end{eqnarray*}
This definition is essential for the techniques considered in this paper.
To illustrate this, we note that the singular value $\sigma_1(K_\mu,\mu)$
associated to a Markov kernel $K$ and a positive probability measure
$\mu$
is the square root of the second largest eigenvalue of
$K^*_\mu K_\mu\dvtx\ell^2(\mu')\ra\ell^2(\mu')$, $\mu'=\mu K$.
This operator is associated with the Markov kernel
\[
P(x,y)=\frac{1}{\mu'(x)}\sum_{z}\mu(z)K(z,x)K(z,y),
\]
which is reversible with respect to $\mu'$ and has associated Dirichlet form
\[
\mathcal E_{P,\mu'}(f,f)=\frac{1}{2}\sum_{x,y,z}|f(x)-f(y)|^2\mu(z)K(z,x)K(z,y).
\]
Hence, using the classical variational formula for eigenvalues, we have
\[
1-\sigma_1(K,\mu)=\inf\biggl\{ \frac{\mathcal E_{P,\mu'}(f,f)}
{\operatorname{Var}_{\mu'}(f)}\dvtx
f\in\ell^2(\mu'), \operatorname{Var}_{\mu'}(f)\neq0\biggr\},
\]
where $\operatorname{Var}_{\mu'}(f)= \|f\|^2_{\ell^2(\mu')}-\mu'(f)^2=
{\sum_x}|f(x)-\mu'(f)|^2\mu'(x)$.
%%%%%%%%%%%%%%%%%%%%%%%%%%%%%%%%%%%%%%%%%%%%%%%%%%%%%%%%%%%%%%%%%%%%%%%%%%%%%%
%%%%%%%%%%%%%%%%%%%%%%%%%%%%%%%%%%%%%%%%%%%%%%%%%%%%%%%%%%%%%%%%%%%%%%%%%%%%%%

%s2.4 ###
\subsection{Nash inequalities}\label{section-Nash}

The use of Nash inequalities to study the convergence of ergodic
(time homogeneous) finite Markov chains was developed in \cite{DS-N}
(Section 7 of \cite{DS-N} discusses time homogeneous chains that admits
an invariant measure).
We refer the reader to that paper for background on this technique.
In this section, we observe that it can be implemented in the context
of time inhomogeneous chains. We start with some basic material.
\begin{defin}\label{def-opnorm}
Let $V$ be a state space equipped with a Markov kernel $K$ and
probability measures
$\mu$ and $\nu$.
If $1\leq p,q\leq\infty$ then
\[
\|K\|_{\ell^p(\mu)\ra\ell^q(\nu)}=\sup_{\|f\|_{\ell^p(\mu)\leq
1}}\bigl\{\|Kf\|_{\ell^q(\nu)}\bigr\}.
\]
If $p$ and $q$ are conjugate exponents, that is, if $1/p+1/q=1$, then
\[
\|f\|_{\ell^p(\mu)}=\sup_{\|g\|_{\ell^q(\mu)\leq1}}\{\langle f,g\rangle
_\mu\}.
\]
\end{defin}

The following proposition is well known in a much more general
context.
\begin{pro}\label{pro-dualnorm} Let $K$ be a Markov kernel.
Let $K_\mu\dvtx\ell^2(\mu K)\ra\ell^2(\mu)$ be the Markov operator on $V$
with adjoint $K^*_\mu\dvtx\ell^2(\mu)\ra\ell^2(\mu K)$ with respect to
the inner product
\[
\langle Kf,g\rangle_\mu=\langle f,K^* g\rangle_{\mu K}.
\]
If $1\leq p,r,s\leq\infty$, $1/p+1/q=1$ and $1/r+1/s=1$ then
\[
\|K\|_{\ell^p(\mu K)\ra\ell^r(\mu)}=\|K^*\|_{\ell^s(\mu)\ra\ell^q(\mu K)}.
\]
\end{pro}

Let now $(K_i)_1^\infty$ be a sequence of Markov kernels on $V$.
Fix a positive probability measure $\mu_0$ and set $\mu_n=\mu_0K_{0,n}$
as usual.
Consider
$K_i\dvtx \ell^2(\mu_i)\ra\ell^2(\mu_{i-1})$, its adjoint
$K^*_i\dvtx\ell^2(\mu_{i-1})\ra\ell^2(\mu_{i})$ and
$P_i=K_i^*K_i\dvtx\ell^2(\mu_{i})\ra\ell^2(\mu_{i})$. The operator $P_i$
is given by the Markov kernel
%
%e2.6 ###
%
\begin{equation}
P_i(x,y)=\frac{1}{\mu_i(x)}\sum_z \mu_{i-1}(z)K_i(z,x)K_i(z,y).
\end{equation}
This kernel is reversible with reversible measure $\mu_i$. We let
\[
\mathcal E_{P_i,\mu_i}(f,f)=\frac{1}{2}\sum_{x,y}|f(x)-f(y)|^2\mu_i(x)P_i(x,y)
\]
be the associated Dirichlet form on $\ell^2(\mu_i)$.
\begin{theo}\label{thm-nash}
Referring to the setup and notation introduced above,
let \mbox{$N\geq1$} and assume that there are constants
$C,D>0$ such that for $1\leq m\leq N$ the following
Nash inequalities hold
%
%e2.7 ###
%
\begin{eqnarray}\label{eqn-nash}
\forall f \dvtx V\ra\mathbb R\qquad \|f\|_{\ell^2(\mu_m)}^{2+1/D} &\leq&
C\biggl(\mathcal{E}_{P_m,\mu_m}(f,f)\nonumber\\[-8pt]\\[-8pt]
&&\hspace*{13.3pt}{} + \frac{1}{N}\|f\|_{\ell^2(\mu
_m)}^2\biggr)\|f\|_{\ell^1(\mu_m)}^{1/D} .\nonumber
\end{eqnarray}
Then, for $0\leq m\leq n\leq N$,
%
%e2.8 ###
%
\begin{equation}\qquad
\max\bigl\{\|K_{m,n}\|_{\ell^2(\mu_n)\ra\ell^{\infty}(\mu_m)},
\|K_{m,n}\|_{\ell^1(\mu_n)\ra\ell^{2}(\mu_m)}\bigr\}
\leq\biggl(\frac{4CB}{n-m+1}\biggr)^D,
\end{equation}
where $B=B(D,N)=(1+1/N)(1+\lceil4D\rceil)$.
\end{theo}
\begin{pf}
Let $(K_i)_0^{\infty}$ be a sequence of Markov kernels on $V$
such that the Nash inequalities
(\ref{eqn-nash}) hold.
Pick a function $f$ such that $\|f\|_{\ell^1(\mu_n)}=1$. For
$1\leq m\leq n\leq N$ define
\[
t_n(n-m)=\|K_{m,n}f\|^2_{\ell^2(\mu_m)}.
\]
Note that for any $n>0$, $(t_n(i))_{i=0}^n$ is nonincreasing.
Indeed, using the contraction property
(\ref{cont}), we have
\begin{eqnarray*}
t_n(i+1)&=&\|K_{n-i-1,n}f\|^2_{\ell^2(\mu_{n-i-1})}=\|
K_{n-i}K_{n-i,n}f\|_{\ell^2(\mu_{n-i-1})}^2\\
&\leq&\|K_{n-i,n}f\|_{\ell^2(\mu_{n-i})}^2=t_n(i).
\end{eqnarray*}
Moreover, note that for any $0\leq i-1\leq n\leq N$
\[
t_n(i)^{1+1/(2D)}\leq C\bigl(t_n(i)-t_n(i+1)+t_n(i)/N\bigr),
\]
where $C$ and $D$ are the constants in (\ref{eqn-nash}).
This follows by applying the Nash inequality to the function $K_{n-i,n}f$.
Corollary 3.1 of \cite{DS-N} then yields that
\[
t_n(i)\leq\biggl(\frac{CB}{i+1}\biggr)^{2D},\qquad 0\leq i\leq n\leq N,
\]
where $B=B(D,N)=(1+1/N)(1+\lceil4D\rceil)$. In particular,
if $0\leq m\leq n\leq N$,
\[
\|K_{m,n}\|_{\ell^1(\mu_n)\ra\ell^2(\mu_m)}\leq\bigl((CB)/(n-m+1)\bigr)^D.
\]
From Proposition \ref{pro-dualnorm} it follows that,
for $0\leq m\leq n\leq N$,
\[
\|K_{m,n}^*\|_{\ell^2(\mu_m)\ra\ell^{\infty}(\mu_n)}\leq\bigl((CB)/(n-m+1)\bigr)^D.
\]
Next we bound $\|K_{m,n}^*\|_{\ell^1(\mu_m)\ra\ell^{\infty
}(\mu_n)}$ for
$0\leq m\leq n\leq N$. Consider the quantity $M(N)$ where
\[
M(N)=\max_{0\leq m\leq n\leq N}
\bigl\{
(n-m+1)^{2D}\|K_{m,n}^*\|_{\ell^1(\mu_m)\ra\ell^{\infty}(\mu_n)}\bigr\}.
\]
Let $l=\lfloor\frac{n-m}{2}\rfloor+m$,
so that $0\leq m\leq l\leq n\leq N$. We have
\begin{eqnarray*}
\|K_{m,n}^*\|_{\ell^1(\mu_m)\ra\ell^{\infty}(\mu_n)}&\leq&
\|K_{m,l}^*\|_{\ell^1(\mu_m)\ra\ell^2(\mu_l)}\|K_{l,n}^*\|_{\ell^2(\mu
_l)\ra\ell^{\infty}(\mu_n)}\\
&\leq&\biggl(\frac{CB}{n-l+1}\biggr)^D\|K_{m,l}^*\|_{\ell^1(\mu_m)\ra\ell
^2(\mu_l)}.
\end{eqnarray*}
Note that for all $0\leq m\leq l\leq N$
%
%e2.9 ###
%
\begin{equation}\label{norm-2-infty}
\|K_{m,l}^*\|_{\ell^1(\mu_m)\ra\ell^2(\mu_l)}\leq
\|K_{m,l}^*\|_{\ell^1(\mu_m)\ra\ell^{\infty}(\mu_l)}^{1/2}
\|K_{m,l}^*\|_{\ell^1(\mu_m)\ra\ell^1(\mu_l)}^{1/2}.
\end{equation}
This follows from the fact that for any function $f$
\[
\|K_{m,l}^*f\|_{\ell^2(\mu_l)}\leq
\|K_{m,l}^*f\|_{\ell^{\infty}(\mu_l)}^{1/2}
\|K_{m,l}^*f\|_{\ell^1(\mu_l)}^{1/2}.
\]
By (\ref{cont}), we have
\begin{eqnarray*}
\|K_{m,n}^*\|_{\ell^1(\mu_m)\ra\ell^{\infty}(\mu_n)}&\leq&
\biggl(\frac{CB}{n-l+1}\biggr)^D\|K_{m,l}^*\|_{\ell^1(\mu_m)\ra\ell
^{\infty}(\mu_l)}^{1/2}\\
&\leq&\biggl(\frac{CB}{(n-l+1)(l-m+1)}\biggr)^DM(N)^{1/2}\\
&\leq&\biggl(\frac{4CB}{(n-m+1)^2}\biggr)^DM(N)^{1/2}.
\end{eqnarray*}
The last inequality follows from the fact that
\[
n-l+1\geq\frac{n-m+1}{2} \quad\mbox{and}\quad l-m+1\geq\frac{n-m+1}{2}.
\]
So we have $M(N)\leq(4CB)^{2D}$ and it follows that for
all $0\leq m\leq n\leq N$
\[
\|K_{m,n}^*\|_{\ell^1(\mu_m)\ra\ell^{\infty}(\mu_n)}\leq
\biggl(\frac{4CB}{n-m+1}\biggr)^{2D}.
\]
By duality, we get that
\[
\|K_{m,n}\|_{\ell^{1}(\mu_n)\ra\ell^{\infty}(\mu_m)}\leq
\biggl(\frac{4CB}{n-m+1}\biggr)^{2D}.
\]
Next, we use the Riesz--Thorin interpolation theorem, see
\cite{SW}, page 179,
which gives us the desired result.
\end{pf}

The next results show how Theorem \ref{thm-nash} together with
the singular value technique of Section \ref{sec-sv} yields merging
results.
\begin{theo}
\label{thm-nash2} Referring to the above setup and notation,
let $N\geq1$ and assume that there are constants
$C,D>0$ such that for $1\leq m\leq N$ the
Nash inequalities
%
%e2.10 ###
%
\begin{eqnarray}\label{eqn-nash2}
\forall f \dvtx V\ra\mathbb R\qquad
\|f\|_{\ell^2(\mu_m)}^{2+1/D}
&\leq&
C\biggl(\mathcal{E}_{P_m,\mu_m}(f,f)+
\frac{1}{N}\|f\|_{\ell^2(\mu_m)}^2\biggr)\nonumber\\[-8pt]\\[-8pt]
&&{} \times\|f\|_{\ell^1(\mu_m)}^{1/D}\nonumber
\end{eqnarray}
hold. Let $\sigma_1(K_m,\mu_{m-1})$ be the second largest singular value
of $K_m\dvtx\ell^2(\mu_m)\ra\ell^2(\mu_{m-1})$, that is, the square
root of
the second largest eigenvalue of $P_m$.
Then, for $n> m$, $N\ge m\ge0$, we have
%
%e2.11 ###
%
\begin{equation}\label{N+S}
d_2(K_{0,n}(x,\cdot),\mu_n)\leq\biggl(\frac{8C(1+\lceil4D\rceil
)}{(m+1)}\biggr)^D
\prod_{m+1}^n \sigma_1(K_i,\mu_{i-1}).
\end{equation}
Moreover, for any $n=2m+u$, $0\le m\le N$, we have
%
%e2.12 ###
%
\begin{equation}\label{sup-N+S}
\max_{x,y}\biggl\{\biggl|\frac{K_{0,n}(x,y)}{\mu_n(y)}-1\biggr|\biggr\}
\le\biggl(\frac{8C(1+\lceil4D \rceil)}{(m+1)}\biggr)^{2D}
\prod_{m+1}^{m+u} \sigma_1(K_i,\mu_{i-1}).\hspace*{-22pt}
\end{equation}
\end{theo}
\begin{pf}
We have
\[
\max_{x\in V}\{d_2(K_{0,n}(x,\cdot),\mu_n)^2\}=
\|K_{0,n}-\mu_n\|^2_{\ell^2(\mu_n)\ra\ell^\infty(\mu_0)},
\]
where $\mu_n$ is understood as the expectation operator $f\mapsto\mu_n (f)$.
Moreover, for any $0\le m\le n$,
\[
K_{0,n}-\mu_n = K_{0,m}(K_{m,n}-\mu_n),
\]
because $K_{0,m}\mu_n f= K_{0,m}\mu_n(f)=\mu_n(f)$.
Hence, for $0\le m\le N$,
\begin{eqnarray*}
d_2(K_{0,n}(x,\cdot),\mu_n)^2&\le&
\|K_{m,n}-\mu_n\|^2_{\ell^2(\mu_n)\ra\ell^2(\mu_m)}
\|K_{0,m}\|^2_{\ell^2(\mu_m)\ra\ell^\infty(\mu_0)}\\
&\le& \Biggl(\prod_{m+1}^n \sigma_1(K_i,\mu_{i-1})^2\Biggr)
\biggl(\frac{4CB}{m+1}\biggr)^{2D}.
\end{eqnarray*}
Using $B= N^{-1}(N+1)(1+\lceil4D \rceil)$, gives (\ref{N+S}).
To obtain the stronger result
(\ref{sup-N+S}), write
\[
\max_{x,y\in V}\biggl\{\biggl|\frac{K_{0,n}(x,y)}{\mu_n(y)}-1\biggr|
\biggr\}=
\|K_{0,n}-\mu_n\|_{\ell^1(\mu_n)\ra\ell^\infty(\mu_0)}
\]
and
\begin{eqnarray*}
&&\|K_{0,n}-\mu_n\|_{\ell^1(\mu_n)\ra\ell^\infty(\mu_0)}\\
&&\qquad\le
\|K_{n-m,n}\|_{\ell^1(\mu_n)\ra\ell^2(\mu_{n-m})}
\times
\|K_{m,n-m}-\mu_{n-m}\|_{\ell^2(\mu_{n-m})\ra\ell^2(\mu_m)}\\
&&\qquad\quad{}\times\|K_{0,m}\|_{\ell^2(\mu_m)\ra\ell^{\infty}(\mu_0)}.
\end{eqnarray*}
The stated bound (\ref{sup-N+S}) follows.
\end{pf}

Just as we did for singular values, let us emphasize that the powerful
looking results stated in this theorem are actually
extremely difficult to apply. Again, the point is that the Dirichlet form
$\mathcal E_{P_m,\mu_m}$, the space $\ell^2(\mu_m)$, and the singular values
$\sigma_1(K_m,\mu_{m-1})$ all involve the unknown sequence of measures
$\mu_n=\mu_0K_{0,n}$, $n=0,\ldots.$ The following subsection gives
similar but more applicable results under additional hypotheses
involving the
notion of $c$-stability.

%s2.5 ###
\subsection{Nash inequality under $c$-stability}

We state two results that parallel Theorems 5.9 and 5.10 of \cite{SZ3}.
\begin{theo}\label{th-nash-stab}
Fix $c\in(1,\infty)$.
Let $(K_i)_1^\infty$ be a sequence of irreducible
Markov kernels on a
finite set $V$.
Assume that $(K_i)_1^\infty$ is $c$-stable with respect to a
positive probability measure $\mu_0$.
For each $i$, set $\mu_0^i=\mu_0K_i$ and let $\sigma(K_i,\mu_0)$ be the
second largest singular value of $K_i=K_{i,\mu_0}$ as an operator from
$\ell^2(\mu^i_0)$ to $\ell^2(\mu_0)$. Let $P^0_i= K_{i,\mu_0}^*K_{i,\mu_0}$.
Let $N\geq1$ and assume that there are constants
$C,D>0$ such that for $1\leq m\leq N$ the
Nash inequalities
%
%e2.13 ###
%
\begin{eqnarray}\label{eqn-nash3}\quad
\forall f \dvtx V\ra\mathbb R\qquad
\|f\|_{\ell^2(\mu^0_m)}^{2+1/D}
&\leq&
C\biggl(\mathcal{E}_{P^0_m,\mu^0_m}(f,f)+
\frac{1}{N}\|f\|_{\ell^2(\mu^0_m)}^2\biggr)\nonumber\\[-8pt]\\[-8pt]
&&{}\times\|f\|_{\ell^1(\mu^0_m)}^{1/D}\nonumber
\end{eqnarray}
holds. Then, for $n> m$, $N\ge m\ge0$, we have
%
%e2.14 ###
%
\begin{eqnarray}\label{N+S3}
d_2(K_{0,n}(x,\cdot),\mu_n)&\leq&\biggl(\frac{8Cc^{2+3/2D}(1+\lceil
4D\rceil)}{(m+1)}\biggr)^D\nonumber\\[-8pt]\\[-8pt]
&&{}\times
\prod_{m+1}^n
\biggl(1-\frac{1-\sigma(K_i,\mu_0)^2}{c^2}\biggr)^{1/2}.\nonumber
\end{eqnarray}
Moreover, for any $n=2m+u$, $0\le m\le N$, we have
\begin{eqnarray*}
\max_{x,y}\biggl\{\biggl|\frac{K_{0,n}(x,y)}{\mu_n(y)}-1\biggr|\biggr\}
&\le&\biggl(\frac{8Cc^{2+3/2D}(1+\lceil4D\rceil)}{(m+1)}\biggr)^{2D}\\
&&{}\times\prod_{m+1}^{m+u} \biggl(1-\frac{1-\sigma(K_i,\mu_0)^2}{c^2}\biggr)^{1/2}.
\end{eqnarray*}
\end{theo}
\begin{pf}
First note that since $\mu_{i-1}/\mu_0\in[1/c,c]$, we
have $\mu_0^i/\mu_i\in[1/c,c]$.
Consider the operator $P_i$ with kernel
\[
P_i(x,y)= \frac{1}{\mu_i(x)}\sum_{z} \mu_{i-1}(z)K_i(z,x)K_i(z,y).
\]
By assumption
\[
\mu_i(x) P_i(x,y)\ge c^{-1}
\mu_0^i(x)\biggl[\frac{1}{\mu^i_0(x)}\sum_{z} \mu
_0(z)K_i(z,x)K_i(z,y)\biggr],
\]
where the term in brackets on the right-hand side is the kernel of
$P^0_i$.
This kernel has second largest eigenvalue $\sigma(K_i,\mu_0)^2$.
A simple eigenvalue comparison argument yields
\[
1-\sigma_1(K_i,\mu_{i-1})^2 \ge\frac{1}{c^2} \bigl(1- \sigma(K_i,\mu_0)^2\bigr).
\]
Further, comparison of measures and Dirichlet form yields the Nash inequality
\begin{eqnarray*}
\forall f \dvtx V\ra\mathbb R\qquad
\|f\|_{\ell^2(\mu_m)}^{2+1/D}
&\leq&
Cc^{2+3/2D}\biggl(\mathcal{E}_{P_m,\mu_m}(f,f)+
\frac{1}{N}\|f\|_{\ell^2(\mu_m)}^2\biggr)\\
&&{}\times\|f\|_{\ell^1(\mu_m)}^{1/D}.
\end{eqnarray*}
Together with Theorem \ref{thm-nash2}, this gives the stated result.
\end{pf}

The next result is based on a stronger hypothesis.
\begin{theo}\label{th-nash-stab2}
Fix $c\in(1,\infty)$. Let $\mathcal Q$ be a family of irreducible aperiodic
Markov kernels on a finite set $V$. Assume that $\mathcal Q$ is $c$-stable
with respect to some positive probability measure $\mu_0$.

Let $(K_i)_1^\infty$ be a sequence of Markov kernels with
$K_i\in\mathcal Q$ for all $i$. Let $\pi_i$ be the invariant measure
of $K_i$.
Let $\tilde{P}_i=K_i^*K_i$ where $K_i\dvtx\ell^2(\pi_i)\ra\ell
^2(\pi_i)$.
Let $\sigma_1(K_i)$ be the
second largest singular value of $K_i$ as an operator on $\ell^2(\pi_i)$.
Let $N\geq1$ and assume that there are constants
$C,D>0$ such that for $1\leq m\leq N$ the
Nash inequalities
%
%e2.15 ###
%
\begin{eqnarray}\label{eqn-nash4}
\forall f \dvtx V\ra\mathbb R\qquad
\|f\|_{\ell^2(\pi_m)}^{2+1/D}
&\leq&
C\biggl(\mathcal{E}_{\tilde{P}_m,\pi_m}(f,f)+
\frac{1}{N}\|f\|_{\ell^2(\pi_m)}^2\biggr)\nonumber\\[-8pt]\\[-8pt]
&&{}\times\|f\|_{\ell^1(\pi_m)}^{1/D}.\nonumber
\end{eqnarray}
Then, for $n> m$, $N\ge m\ge0$, we have
%
%e2.16 ###
%
\begin{eqnarray}\label{N+S4}
d_2(K_{0,n}(x,\cdot),\mu_n)&\leq&\biggl(\frac{8Cc^{4+3/D}(1+\lceil
4D\rceil)}{(m+1)}\biggr)^D\nonumber\\[-8pt]\\[-8pt]
&&{}\times\prod_{m+1}^n
\biggl(1-\frac{1-\sigma_1(K_i)^2}{c^4}\biggr)^{1/2}.\nonumber
\end{eqnarray}
Moreover, for any $n=2m+u$, $0\le m\le N$, we have
\[
\max_{x,y}\biggl\{\biggl|\frac{K_{0,n}(x,y)}{\mu_n(y)}-1\biggr|\biggr\}
\le\biggl(\frac{8Cc^{4+3/D}(1+\lceil4D\rceil)}{(m+1)}\biggr)^{2D}
\prod_{m+1}^{m+u} \biggl(1-\frac{1-\sigma_1(K_i)^2}{c^4}\biggr)^{1/2}.
\]
\end{theo}
\begin{pf} Note that the hypothesis that $\mathcal Q$ is $c$-stable
implies $\pi_i/\mu_j\in[1/c^2,c^2]$ for all $i,j$.
Consider again the operator $P_i$ and its kernel
\[
P_i(x,y)= \frac{1}{\mu_i(x)}\sum_{z} \mu_{i-1}(z)K_i(z,x)K_i(z,y).
\]
By assumption
\begin{eqnarray*}
\mu_i(x) P_i(x,y) &\ge& c^{-2}
\pi_i(x)\biggl[\frac{1}{\pi_i(x)}\sum_{z} \pi
_{i}(z)K_i(z,x)K_i(z,y)\biggr]\\
&\ge& c^{-2} \pi_i(x)\tilde{P}_i(x,y).
\end{eqnarray*}
A comparison argument
similar to the one used in the previous proof yields the desired result.
\end{pf}

%%%%%%%%%%%%%%%%%%%%%%%%%%%%%%%%%%%%%%%%%%%%%%%%%%%%%%%%%%%%%%%%%%%%%%%%%%%%

%s3 ###
\section{Examples involving Nash inequalities}
This section describes applications of the Nash inequality technique to
several examples. All these examples are of the following general type.

(1) There is a basic reversible model $(K,\pi)$ on a space $V_N$
(growing with $N$) that is well understood
because:
\begin{itemize}
\item We have good grasp on the second largest
singular value $\sigma_N$ of $(K,\pi)$.
\item The model $(K,\pi)$ satisfies a good
Nash inequality, that is, an inequality of the form
\[
\|f\|^{2+1/D}_{\ell^2(\pi)}\le
BT_N \biggl(\mathcal E_{K^*K,\pi}(f,f)+\frac{1}{bT_N}
\|f\|_{\ell^2(\pi)}^2\biggr)\|f\|_{\ell^1(\pi)}^{1/D}
\]
with $B,b$ independent of $N$ and $T_N\simeq(1-\sigma_N)^{-1}$. Here,
$f\simeq g$ implies that there exist constants $d,D>0$ such that
$dg\leq f\leq Dg$.
\item Together, the Nash inequality and second largest singular value estimate
yield the mixing time estimate
\[
\max_{x,y}\biggl\{\biggl|\frac{K^t(x,y)}{\pi(y)}-1\biggr|\biggr\}\le\eta,\qquad
t\ge\frac{A(1+\log_+ 1/\eta)}{1-\sigma_N},
\]
where $A$ is independent of $N$.
\end{itemize}

(2) We are given a sequence $(K_i)_1^\infty$ or a set
$\mathcal Q_N$ of Markov kernels on $V_N$ which satisfies:
\begin{itemize}
\item$(K_i)_1^\infty$ or $\mathcal Q_N$ is $c$-stable
with respect to a measure $\mu_0$ which is either equal or at least
comparable to $\pi$.
\item The Markov kernels $K_i$ or the elements of $\mathcal Q_N$ are all
bounded perturbations of $K$ in the sense that
$K_i(x,y)/K(x,y)$ is bounded away from $0$ and away from $\infty$
for all $(x,y)\in V_N^2$. In particular, $K_i(x,y)=0$ if and only if
$K(x,y)=0$.
\end{itemize}

Under such circumstances,
Theorem \ref{th-nash-stab} (or Theorem \ref{th-nash-stab2}) applies and
yields the conclusion that the time inhomogeneous Markov chain associated
with the sequence $K_i$ under investigation has a relative-sup
merging time $T_\infty(\eta)$ bounded by
\[
T_\infty(\eta)\le\frac{A'(1+\log_+ 1/\eta)}{1-\sigma_N}
\]
for some constant $A'$ independent of $N$.

The most obvious basic model is, perhaps, the simple random walk
on $\mathbb Z/ N\mathbb Z$ (with some holding
if $N$ is even to avoid periodicity). This model has $1-\sigma_N\simeq1/N^2$
and satisfies the desired Nash inequality with $D=1/4$. The first
subsection presents applications to a perturbation of this model.

%s3.1 ###
\subsection{Asymmetric perturbation at the middle vertex}

\begin{figure}%[b]

\includegraphics{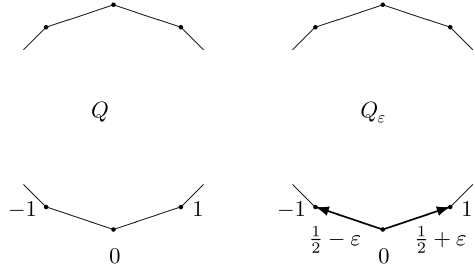}

\caption{The asymmetric perturbation.}\label{fig1}
%
%
%
%
%{\thicklines\put(150,0){\vector(-3,1){30}}
%
\end{figure}

In this example, $V_N=\mathbb Z/ p_N\mathbb Z$ is a finite circle.
It will be convenient to enumerate the points in $V_N$ by writing $
V_N=\{-(N-1),\ldots,-1,0,1,\ldots,(N-1),N\}$ if $p_N=2N$ and
$V_N=\{-N,\ldots,-1,0,1,\ldots,N\}$ if $p_N=2N+1$.
The simple random walk in $V$ has kernel
%
%e3.1 ###
%
\begin{equation}\label{SRW-circ}
Q(x,y)=
\cases{
1/2, &\quad if $|x-y|=1$, \cr
0, &\quad otherwise,}
\end{equation}
and reversible measure $u\equiv\frac{1}{p_N}$.
For any $\varepsilon>0$, define the perturbation kernel
%
%e3.2 ###
%
\begin{equation}
\Delta_{\varepsilon}(x,y)=
\cases{
\varepsilon, &\quad if $(x,y)=(0,1)$, \cr
-\varepsilon, &\quad if $(x,y)=(0,-1)$, \cr
0, &\quad otherwise.}
\end{equation}
For $\varepsilon\in(-1/2,1/2)$, the Markov kernel
$Q_{\varepsilon}=Q+\Delta_{\varepsilon}$
is a perturbation of $Q$. See Figure \ref{fig1}.

%f1 ###
%

For any fixed $0<\varepsilon<1/2$, set
\[
\mathcal{Q}(\varepsilon)=\{Q_{\delta}\dvtx \delta\in[-\varepsilon
,\varepsilon]\}.
\]
We shall see below that $\mathcal{Q}(\varepsilon)$ is $c$-stable.
\begin{defin}
Let $\mathcal S_N(\varepsilon)$ be the set of all probability
measures on $V_N$ which satisfy the following two properties:
\begin{enumerate}[(2)]
\item[(1)] for all $x\in V_N$, there exist constants $a_{\mu,x}$ such
that $a_{\mu,x}=-a_{\mu,-x}$ and
\[
\mu(x)=(1/p_N)+a_{\mu,x}
\]
\item[(2)] for all $x\in V_N$ we have that $|a_{\mu,x}|\leq2\varepsilon/p_N$.
\end{enumerate}
\end{defin}
\begin{rem}
Note that we always have $a_{\mu,0}=0$ (since $-0=0$) and,
in the case when $p_N=2N$, $a_{\mu,N}=0$.
\end{rem}
\begin{claim}\label{thm-cstab}
Let $\mu\in\mathcal S_N(\varepsilon)$ defined above,
then for any $K \in\mathcal Q(\varepsilon)$ we have
that $\mu K \in\mathcal S_N(\varepsilon)$.
\end{claim}
\begin{pf}
Let $\mu\in\mathcal S_N(\varepsilon)$ and $K=Q_{\delta}\in\mathcal
Q(\varepsilon)$,
$\delta\in[-\varepsilon,\varepsilon]$. We show that $\mu K$ has the
properties required to be in $\mathcal S_N(\varepsilon)$.

(1) Any measure $\mu\in\mathcal S_N$ can be written as
$\mu=u+m_{\mu}$
where $m_{\mu}$ is the (nonprobability) measure $m_{\mu}(x)=a_{\mu,x}$.
A simple calculation yields that
\[
m_{\mu}Q(x)=(a_{\mu,x-1}+a_{\mu,x+1})/2.
\]
Since $a_{\mu,x}=-a_{\mu,-x}$, we obtain that
\[
m_{\mu}Q(x)=-m_{\mu}Q(-x) \quad\mbox{and}\quad m_{\mu}Q(0)=0.
\]
The fact that $\mu Q=(u+m_\mu)Q=u + m_\mu Q$ implies
that $\mu Q$ satisfies property $(1)$ in the definition of $\mathcal
S_N(\varepsilon)$.
To see that $\mu Q_{\delta}\in\mathcal S_N(\varepsilon)$ also satisfies
this property,
we note that
\[
\mu\Delta_{\delta}(x)=
\cases{
\delta\mu(0), &\quad if $x=1$, \cr
-\delta\mu(0), &\quad if $x=-1$, \cr
0, &\quad otherwise.}
\]
It now follows that $\mu Q_{\delta}\in\mathcal S_N$ has property $(1)$
in the
definition of $\mathcal S_N(\varepsilon)$ since $\mu Q_{\delta}=
\mu(Q+\Delta_{\delta})$.

(2) We consider the measure $\mu K$.
For $x\notin\{-1,1\}$ property $(2)$ of $\mathcal S_N(\varepsilon)$
follows easily from the fact that $|a_{\mu,x}|\leq2\varepsilon/p_N$ and
\[
\mu K(x)=1/p_N+1/2(a_{\mu,{x-1}}+a_{\mu,{x+1}}).
\]
For $x=1$, we note that
\begin{eqnarray*}
\mu K(1)&\leq&\mu(0)(1/2+\varepsilon)+\mu(2)(1/2) = 1/p_N+\varepsilon
/p_N+(1/2)a_{\mu,2}\\
&\leq& 1/p_N+2\varepsilon/p_N.
\end{eqnarray*}
Similarly
\begin{eqnarray*}
\mu K(1)&\geq&\mu(0)(1/2-\varepsilon)+\mu(2)(1/2) = 1/p_N-\varepsilon
/p_N-(1/2)a_{\mu,2}\\
&\geq& 1/p_N-2\varepsilon/p_N.
\end{eqnarray*}
The proof now follows from the fact that $a_{\mu K,1}=-a_{\mu K,-1}$ as
proved in
part (1) above.
\end{pf}
\begin{claim}
The family $\mathcal Q(\varepsilon)$ is $\frac{1+2\varepsilon
}{1-2\varepsilon}$-stable
with respect to any $\mu_0\in\mathcal S_N(\varepsilon)$.
\end{claim}
\begin{pf}
Claim \ref{thm-cstab} implies that for any sequence $(K_i)_0^{\infty}$ such
that $K_i\in\mathcal{Q}_{\varepsilon}$ and any measure
$\mu_0\in\mathcal S_N(\varepsilon)$
we have $\mu_n=\mu_0K_{0,n}\in\mathcal S_N(\varepsilon)$ for all $n\geq0$.
Note that for any measure $\nu\in\mathcal S_N(\varepsilon)$ we have that
\[
\nu(x)=1/p_N+a_{\nu,x}\leq(1+2\varepsilon)/p_N \quad\mbox{and}\quad
\nu(x)=1/p_N+a_{\nu,x}\geq(1-2\varepsilon)/p_N.
\]
Hence,
\[
\frac{1-2\varepsilon}{1+2\varepsilon}\leq
\frac{\mu_n(x)}{\mu_0(x)}
\leq\frac{1+2\varepsilon}{1-2\varepsilon}.
\]
\upqed\end{pf}

When
$p_N=2N$, the kernels $Q_\delta$ yield periodic chains on $V_N$.
In this case, we will study the merging properties of
\[
\mathcal Q_{\mathrm{lazy}}
(\varepsilon)= \bigl\{ \tfrac{1}{2}(I+K)\dvtx K\in\mathcal Q(\varepsilon)
\bigr\},
\]
that is, the so-called lazy version of $\mathcal Q(\varepsilon)$.
We set
\[
\overline{Q}_\delta= \tfrac{1}{2}(I+Q_\delta).
\]
For any
$\mu\in\mathcal S_N(\varepsilon)$, we consider the
kernel
\[
P_{\delta,\mu}(x,y)=
\frac{1}{\mu\overline{Q}_\delta(x)}\sum_z \mu(z)\overline{Q}_\delta(z,x)
\overline{Q}_\delta(z,y),
\]
which is the kernel of $K^*K$ where $K=\overline{Q}_\delta\dvtx
\ell^2(\mu\overline Q_\delta)\ra\ell^2(\mu)$. This is $0$ unless
$y=x, x\pm1,x\pm2$
and we compare it to
\begin{eqnarray*}
P(x,y)&=&P_{0,u}(x,y)=
\frac{1}{u(x)}\sum_z u(z)\overline{Q} (z,x)\overline{Q}(z,y)\\
&=& \sum_z \overline{Q}(z,x)\overline{Q}(z,y),
\end{eqnarray*}
which is $3/8$ if $y=x$, $1/4$ if $y=x\pm1$, $1/16$ if $y=x\pm2$
and $0$ otherwise.
The definitions of $Q_\delta$ and $\mathcal S_N(\varepsilon)$ yield
\begin{eqnarray*}
\mu\overline{Q}_{\delta} (x)P_{\delta,\mu}(x,y) &\ge&
\frac{(1-2\varepsilon)(1-2\delta)^2}{(1+2\varepsilon)} u(x)P(x,y)\\
&\ge&
\frac{(1-2\varepsilon)^3}{(1+2\varepsilon)} u(x) P(x,y).
\end{eqnarray*}
This yields
%
%e3.3 ###
%
\begin{equation}\label{eqn-lazycirc-comp-Dir}
\mathcal E_{P,\mu}(f,f)\le\frac{(1+2\varepsilon)}{(1-2\varepsilon)^3}
\mathcal E_{P_{\delta,\mu},\mu\overline Q_{\delta}}(f,f),
\end{equation}
whereas the stability property implies that the relevant measures
$\mu\overline{Q}_\delta$ and $u$ satisfy
%
%e3.4 ###
%
\begin{equation}\label{eqn-lazycirc-comp-meas}
\frac{(1-2\varepsilon)}{(1+2\varepsilon)} u \le
\mu\overline{Q}_\delta\le\frac{(1+2\varepsilon)}{(1-2\varepsilon)} u.
\end{equation}

In the case when $p_N=2N+1$, we may work directly with
the kernels $Q_{\delta}$ as they are not periodic.
An analysis similar to that above will
give versions of (\ref{eqn-lazycirc-comp-Dir}) and
(\ref{eqn-lazycirc-comp-meas}) for $Q_\delta$.

Applying the line of reasoning explained at the beginning of this section
and using Theorem \ref{th-nash-stab2}, we get the following result.
\begin{theo}\label{thm-lazycirc} Fix $\varepsilon\in(0,1/2)$.
For any $\eta>0$ the total variation $\eta$-merging time of the
family $\mathcal Q_{\mathrm{lazy}}(\varepsilon)$ on
$V_N=\mathbb{Z}/2N\mathbb{Z}$
[resp., $\mathcal Q(\varepsilon)$ on $V_N=\mathbb{Z}/(2N+1)\mathbb{Z}$]
is at most
$ B(\varepsilon) N^2( 1+\log_+ 1/\eta)$ for some constant
$B(\varepsilon)\in(0,\infty)$. In fact, we can choose $B(\varepsilon)$
such that
\[
\forall n\ge B(\varepsilon)N^2(1+\log_+ 1/\eta)\qquad
\max_{x,y\in V_N}\biggl\{
\biggl|\frac{K_{0,n}(x,z)}{K_{0,n}(y,z)}-1\biggr|\biggr\}
\le\eta
\]
for any sequence $K_i\in\mathcal Q_{\mathrm{lazy}}(\varepsilon)$
[resp., $K_i\in\mathcal Q(\varepsilon)$].
\end{theo}

%%%%%%%%%%%%%%%%%%%%%%%%%%%%%%%%%%%%%%%%%%%%%%%%%%%%%%%%%%%%%%%%

%s3.2 ###
\subsection{Perturbations of some birth and death chains}
In \cite{SC-N}, Nash inequalities are used to study certain
birth and death chains on $V_N=\{-N,\ldots,0,\break\ldots,N\}$
with reversible measures which belong to one of the following two families:
\[
\hat{\pi}_\alpha(x)= \hat{c}(\alpha,N) (N-|x|+1)^\alpha,\qquad \alpha\ge0,
\]
and
\[
\check{\pi}_\alpha(x)=\check{c}(\alpha,N) (|x|+1)^\alpha,\qquad \alpha\ge0.
\]
Here, we consider $\alpha\in[0,\infty)$ to be a fixed parameter and are
interested in what happens when $N$ tends to infinity. From this perspective,
the normalizing constants $\hat{c}(\alpha,N), \check{c}(\alpha,N)$
are comparable and behave as
\[
\hat{c}(\alpha,N)\simeq\check{c}(\alpha,N)\simeq N^{-\alpha-1}.
\]
Set
\[
\zeta(\alpha,N)=\sum_0^N (1+i)^{-\alpha}\simeq
\cases{
1, &\quad if $\alpha>1$, \cr
\log N, &\quad if $\alpha=1$, \cr
N^{-\alpha+1}, &\quad if $\alpha\in[0,1)$.}
\]

Here, all $\simeq$ must be understood for fixed $\alpha$ and the implied
comparison constants depend on $\alpha$.
Let $\hat{M}_\alpha$ (resp., $\check{M}_\alpha$) be the Markov kernel of
the Metropolis chain with basis the symmetric simple random walk on $V_N$
with holding $1/3$ at all points except at the end points where the
holding is
$2/3$, and target $\hat{\pi}_\alpha$,
(resp., $\check{\pi}_\alpha$). Let $\hat{\lambda}(\alpha,N)$,
$\check{\lambda}(\alpha,N)$ be the corresponding spectral gaps.
Let $\hat{T}(\alpha,N,\eta)$, $\check{T}(\alpha,N,\eta)$
be the relative-sup mixing times of these chains. It is proved in \cite{SC-N}
that
\[
\hat{\lambda}(\alpha,N) \simeq1/N^2,\qquad
\hat{T}(\alpha,N,\eta)\simeq N^2 (1+\log_+ 1/\eta),
\]
whereas
\begin{eqnarray*}
\check{\lambda}(\alpha,N) &\simeq&\check{c}(\alpha,N)/\zeta(\alpha,N) ,\\
\check{T}(\alpha,N,\eta)&\simeq&\bigl(N^2+ [\check{c}(\alpha,N)/\zeta(\alpha,N)]
\log_+ 1/\eta\bigr) .
\end{eqnarray*}
Note that
\[
\check{c}(\alpha,N)/\zeta(\alpha,N)
\simeq
\cases{
N^{-(1+\alpha)}, &\quad if $\alpha>1$, \cr
(N^2 \log N)^{-1}, &\quad if $\alpha=1$, \cr
N^{-2}, &\quad if $\alpha\in[0,1)$.}
\]
These results are based on the Nash inequalities satisfied by these
chains. Namely, letting
$\mathcal E_\alpha=\mathcal E_{\hat{M}_\alpha,\hat{\pi}_\alpha}$
or $\mathcal E_\alpha=\mathcal E_{\check{M}_\alpha,\check{\pi}_\alpha}$
and $\pi_\alpha=\hat{\pi}_\alpha$ or $\pi_\alpha=\check{\pi}_\alpha$,
there are constants $A_\alpha,a_\alpha\in(0,\infty)$ such that
\[
\|f\|_{\ell^2(\pi_{\alpha})}^{2+1/D_\alpha}\le A_\alpha
N^2 \biggl(\mathcal E_{\alpha}(f,f)+\frac{1}{a_\alpha N^2}\|f\|^2_{\ell^2(
\pi_{\alpha})}
\biggr)\|f\|_{\ell^1(\pi_{\alpha})}^{1/D_\alpha}
\]
with $D_\alpha=1+\alpha$. See \cite{SC-N}.

In cite \cite{SZ3}, the authors consider the class of birth and death chains
$Q$ on $V_N=\{-N,\ldots,0,\ldots,N\}$ that are symmetric with respect to
the middle point,
that is, satisfy $Q(x,x+1)=Q(-x,-x-1)$, $Q(x,x-1)=Q(-x,-x+1)$,
$Q(x,x)=Q(-x,-x)$,
$x\in\{0,N\}$. For any such chain $Q$, let $\nu$ be the reversible measure.
It satisfies $\nu(x)=\nu(-x)$.
Consider the perturbation set
\[
\mathcal Q_N(Q,\varepsilon)=\{ Q+\Delta_s\dvtx s\in[-\varepsilon
,\varepsilon]\},\qquad
\varepsilon\in[0,q_0),
\]
where $q_0=Q(0,\pm1)$, $\Delta_s(0,\pm1)=\pm s$
and $\Delta(x,y)=0$ otherwise. These perturbations at the
middle vertex have reversible measure $\nu_s$ that satisfy
\[
\nu_s(0)=\nu(0),\qquad \nu_s(\pm x) = \nu(\pm x)(1 \pm s/q_0),\qquad
x\in\{1,\ldots,N\}.
\]
The main point of this construction is the following.
\begin{pro} Fix $Q$, $\nu$ as above and $\varepsilon\in[0,q_0)$.
The set $\mathcal Q_N(Q,\varepsilon)$ is $c$-stable with respect to
$\mu_0=\nu$ with $c=(q_0+\varepsilon)/(q_0-\varepsilon)$.
\end{pro}

In order to apply this results to our example
$\hat{M}_\alpha,\check{M}_\alpha$, we observe that
\[
\hat{q}_0(\alpha)=\hat{M}_\alpha(0,-1)=
\frac{1}{3}\biggl(\frac{N}{N+1}\biggr)^\alpha
\]
and
\[
\check{q}_0(\alpha)=\check{M}_\alpha(0,-1)= \tfrac{1}{3}.
\]
Now, Theorem \ref{th-nash-stab2} yields the following result.
\begin{theo} Fix $\alpha\in[0,\infty)$ and
set $\hat{\varepsilon}_{N,\alpha} =\frac{1}{6}(N/(N+1))^\alpha$,
$\check{\varepsilon}_{N,\alpha}=1/6$.
\begin{enumerate}
\item There exists a constant $A$ independent of $N$ such that, for
any sequence $(K_i)_1^\infty$ with $K_i\in
\mathcal Q_N(\hat{M}_\alpha,\hat{\varepsilon}_{N,\alpha})$, we have
\[
T_\infty(\eta)\le AN^2 (1+\log_+ 1/\eta).
\]
\item There exists a constant $A$ independent of $N$ such that, for
any sequence $(K_i)_1^\infty$ with $K_i\in
\mathcal Q_N(\check{M}_\alpha,\check{\varepsilon}_{N,\alpha})$, we have
\[
T_\infty(\eta)\le A
\cases{
N^2+ N^{1+\alpha}\log_+ 1/\eta, &\quad if $\alpha>1$, \cr
N^2+ (N^2\log N)\log_+ 1/\eta, &\quad if $\alpha=1$, \cr
N^2(1+\log_+ 1/\eta), &\quad if $\alpha\in(0,1)$.}
\]
\end{enumerate}
\end{theo}

%%%%%%%%%%%%%%%%%%%%%%%%%%%%%%%%%%%%%%%%%%%%%%%%%%%%%%%%%%%%%%%%%%%%%%%%%%%%%%%%%%%%%%%%%%%

%s4 ###
\section{Logarithmic Sobolev inequalities}\label{sec-LS}
This section develops the technique of logarithmic Sobolev inequality
for time inhomogeneous finite Markov chains. It should be noted that
the logarithmic Sobolev technique has been mostly applied in the literature
in the context of continuous time chains. In \cite{Mic},
Miclo tackled the problem of adapting this technique to discrete time
(time homogeneous) chains. There are two different ways to use
logarithmic Sobolev inequality for mixing estimates. One, the most
powerful, provides results for relative-sup merging and is based on
hypercontractivity. The other is based on entropy and only produces
bounds for total variation merging.
We will discuss and illustrate both approaches below in the context
of time inhomogeneous chains. The entropy approach is
already treated in \cite{DLM}.

%%%%%%%%%%%%%%%%%%%%%%%%%%%%%%%%%%%%%%%%%%%%%%%%%%%%%%%%%%%%%%%%%%%%%%%%%%%%%%%%%%%%%%%%%%%

%s4.1 ###
\subsection{Hypercontractivity}\label{subsec-Hyper}

Recall that, for any positive probability distribution~$\mu$,
a Markov kernel $K$ can be thought of as a contraction
\[
K_\mu\dvtx\ell^2(\mu')\ra\ell^2(\mu)\qquad \mbox{for $\mu'=\mu K$.}
\]
The adjoint $K_\mu^*\dvtx\ell^2(\mu)\ra\ell^2(\mu')$ has kernel
\[
K^*_{\mu}(x,y)=\frac{K(y,x)\mu(y)}{\mu'(x)}.
\]
Set $P=K_\mu^*K_\mu\dvtx\ell^2(\mu')\ra\ell^2(\mu')$.
We define the logarithmic Sobolev constant
\[
l(P)=\inf
\biggl\{\frac{\mathcal{E}_{P,\mu'}(f,f)}{\mathcal L(f^2,\mu')}\dvtx
\mathcal L(f^2,\mu')
\neq0, f\neq\mbox{constant}\biggr\},
\]
where the $\ell^2$ relative entropy $\mathcal L(f^2,\nu)$
of a function $f$ with respect to the measure
$\nu$ is defined by
\[
\mathcal L(f^2,\nu)=\sum_{x\in V}
f^2\log\biggl(\frac{f^2}{\|f\|^2_{\ell^2(\nu)}}\biggr)\nu(x).
\]

The following proposition is a slight generalization of
\cite{Mic}, Proposition 2, in that it allows for the necessary
change of measure.
\begin{pro}\label{pro-opnorm-LS}
Let $K$ and $\mu$ be a Markov kernel and a probability measure, respectively.
For all $q_0\geq2$ and $q\leq[1+l(P)]q_0$,
then
\[
\|K\|_{\ell^{q_0}(\mu')\ra\ell^q(\mu)}\leq1.
\]
\end{pro}

In order to prove the proposition above, we will need the following
two lemmas from \cite{Mic}.
\begin{lem}[(\cite{Mic}, Lemma 3)]\label{lem3}
Let $\nu$ be a probability measure.
For all $q\geq\break q_0\geq1$,
\[
\|f\|_{\ell^q(\nu)}-\|f\|_{\ell^{q_0}(\nu)}\leq
\frac{q-q_0}{q_0q}\|f\|_{\ell^q(\nu)}^{1-q}\mathcal L(f^{q/2},\nu).
\]
\end{lem}
\begin{lem}[(\cite{Mic}, Lemma 4)]\label{lem4}
Fix $\nu\geq0$ and $q\geq2$, then for any $t\geq0$ and
$-t\leq s\leq\nu t$ we have that
\[
(t+s)^q\geq t^q+qt^{q-1}s+g(q,\nu)\bigl((t+s)^{q/2}-t^{q/2}\bigr)^2,
\]
where
\[
g(q,v)=\frac{(1+\nu)^q-1-q\nu}{((1+\nu)^{q/2}-1)^2}.
\]
\end{lem}

The proof of Proposition \ref{pro-opnorm-LS} follows directly that
of Proposition 2 in \cite{Mic}.
\begin{pf*}{Proof of Proposition \ref{pro-opnorm-LS}}
To prove Proposition \ref{pro-opnorm-LS} is suffices to only
consider positive functions. For $f> 0$, we begin by writing
%
%e4.1 ###
%
\begin{eqnarray}\label{eqn1-LS}
\|K f\|_{\ell^q(\mu)}-\|f\|_{\ell^2(\mu')}&=&
\|Kf\|_{\ell^q(\mu)}-\|f\|_{\ell^q(\mu')}\nonumber\\[-8pt]\\[-8pt]
&&{}+\|f\|_{\ell^q(\mu')}-\|f\|
_{\ell^2(\mu')}.\nonumber
\end{eqnarray}
The difference of the last two terms on the right-hand side is
controlled by Lem\-ma~\ref{lem3}. To control the first two terms,
we will use the concavity result
\[
\forall a,b\geq0\qquad a^{1/q}-b^{1/q}\leq\frac{1}{q} b^{1/q-1}(a-b).
\]
It follows that
\[
\|Kf\|_{\ell^q(\mu)}-\|f\|_{\ell^q(\mu')}\leq\frac{1}{q}
\|f\|_{\ell^q(\mu')}^{1-q}\bigl(\|K f\|_{\ell^q(\mu)}^q-\|f\|_{\ell
^q(\mu')}^q\bigr).
\]
Set
\[
\nu(K)=\max\{1/K(x,y) \dvtx K(x,y)>0\}-1.
\]
Following the notation of Lemma \ref{lem4}, fix $x,y\in V$ and set
$\nu=\nu(K)$, $t=Kf(x)$ and $t+s=f(y)$. If $K(x,y)>0$, then $-t\leq
s\leq\nu t$ and so
\begin{eqnarray*}
f(y)^q &\geq& Kf(x)^q+qKf(x)^{q-1}\bigl(f(y)-Kf(x)\bigr)\\
&&{} + g(q,\nu(K))
\bigl(f(y)^{q/2}-Kf(x)^{q/2}\bigr)^2.
\end{eqnarray*}
Fix $x$ and integrate with respect to the measure $K(x,\cdot)$ to get
\[
Kf^q(x)\geq(Kf(x))^q+g(q,\nu(K))
\sum_{y\in V}K(x,y)\bigl(f(y)^{q/2}-Kf(x)^{q/2}\bigr)^2.
\]
We also have
\begin{eqnarray*}
\sum_{y\in V}K(x,y)\bigl(f^{q/2}(y)-(Kf(x))^{q/2}\bigr)^2&\geq&
\min_{c\in\mathbb{R}}\sum_{y\in V}K(x,y)\bigl(f^{q/2}(y)-c\bigr)^2\\
&=&\sum_{y\in V}K(x,y)\bigl(f^{q/2}(y)-K(f^{q/2})(x)\bigr)^2\\
&=&Kf^q(x)-(Kf^{q/2}(x))^2.
\end{eqnarray*}
Hence,
\[
Kf^q(x)\geq(Kf(x))^q+g(q,\nu(K))\bigl(Kf^q(x)-(Kf^{q/2}(x))^2\bigr).
\]
Integrating with respect to $\mu$ gives us that
%
%e4.2 ###
%
\begin{equation}\label{eqn2-LS}
\|f\|_{\ell^q(\mu')}^q
\geq\|Kf\|_{\ell^q(\mu)}^q+g(q,\nu(K))
\mathcal E_{P,\mu'}(f^{q/2},f^{q/2}).
\end{equation}
It follows from Lemma \ref{lem3}, (\ref{eqn1-LS}) and (\ref{eqn2-LS}) that
\begin{eqnarray*}
&&\|Kf\|_{\ell^q(\mu)}-\|f\|_{\ell^2(\mu')}\\
&&\qquad\leq\frac{1}{q} \|f\|_{\ell^q(\mu')}^{1-q}
\biggl(\frac{q-q_0}{q_0}\mathcal L(f^{q/2},\mu')
-g(q,\nu(K))\mathcal E_{P,\mu'}(f^{q/2},f^{q/2})\biggr).
\end{eqnarray*}
In \cite{Mic}, it is noted that for all $\nu>0$ and $q\geq2$ we have
$g(q,\nu)\geq1$.
So if $q\leq[1+l(P)]q_0$ then $q\leq[1+g(q,\nu(K))l(P)]q_0$. Hence,
\begin{eqnarray*}
&&\|Kf\|_{\ell^q(\mu)}-\|f\|_{\ell^2(\mu')}\\
&&\qquad\leq\frac{1}{q} \|f\|_{\ell^q(\mu')}^{1-q}g(q,\nu(K))
\bigl(l(P)\mathcal
L(f^{q/2},\mu')-\mathcal E_{P,\mu'}(f^{q/2},f^{q/2})
\bigr).
\end{eqnarray*}
Since $l(P)$ is the logarithmic Sobolev constant, we
get our desired result,
\[
\|Kf\|_{\ell^q(\mu)}-\|f\|_{\ell^2(\mu')}\leq0.
\]
\upqed\end{pf*}
\begin{cor}\label{cor-Sob}
Let $(K_n)_0^{\infty}$ be a sequence of Markov kernels on a finite set
$V$ and
$\mu_0$ be an initial distribution on $V$. Set $\mu_n=\mu_0K_{0,n}$.
Consider $K_i\dvtx\ell^2(\mu_i)\ra\ell^2(\mu_{i-1})$ and
$P_i=K_i^*K_i\dvtx\ell^2(\mu_i)\ra\ell^2(\mu_i)$. Let $l(P_i)$ be the
logarithmic
Sobolev constant of $P_i$.
Then for any $q_0\geq2$ and
$q\leq\prod_{i=1}^n(1+l(P_i))q_0$, we have that
\[
\|K_{0,n}\|_{\ell^{q_0}(\mu_n)\ra\ell^q(\mu_0)}\leq1.
\]
\end{cor}
\begin{pf}
When $n=2$, set $q_1=(1+l(P_2))q_0$, then $q=(1+l(P_1))q_1$. It follows from
Proposition \ref{pro-opnorm-LS} that
\[
\|K_{0,2}\|_{\ell^{q_0}(\mu_2)\ra\ell^q(\mu_0)}\leq
\|K_2\|_{\ell^{q_0}(\mu_2)\ra\ell^{q_1}(\mu_1)}\|K_1\|_{\ell^{q_1}(\mu
_1)\ra\ell^{q_2}(\mu_0)}\leq
1.
\]
The proof by induction follows similarly.
\end{pf}

We now relate the results above to bounds on merging times.
\begin{theo}\label{thm-LS}
Let $V$ be a finite set equipped with a sequence of Markov kernels
$(K_n)_0^{\infty}$
and an initial distribution $\mu_0$. Let $\mu_n=\mu_0K_{0,n}$.
Consider $K_i\dvtx\break\ell^2(\mu_i)\ra\ell^2(\mu_{i-1})$ and
$P_i=K_i^*K_i\dvtx\ell^2(\mu_i)\ra\ell^2(\mu_i)$.
Let $l(P_i)$ be the logarithmic Sobolev constant of $P_i$. Set
\[
m_x=\min\Biggl\{t\in\mathbb{N}\dvtx\sum_{i=1}^t\log\bigl(1+l(P_i)\bigr)\geq\log\log(\mu
_0(x)^{-1/2})\Biggr\}.
\]
Then for $n\geq m_x$, we have that
\[
d_2(K_{0,n}(x,\cdot),\mu_n)^2\leq e^{2}\prod_{i=m_x+1}^n\sigma_1(K_i,\mu
_{i-1})^2.
\]
\end{theo}
\begin{pf}
Fix $x$, and let $m=m_x$. If $0\leq m\leq n$, $K_{0,n}^*=K_{m,n}^*K_{0,m}^*$.
Indeed, for any $f\in\ell^2(\mu_0)$ and $g\in\ell^2(\mu_n)$ we have that
\[
\langle K_{0,n}^*f,g\rangle_{\mu_n}=\langle f,K_{0,n}g\rangle_{\mu_0}=
\langle K_{0,m}^*f,K_{m,n}g\rangle_{\mu_m}=\langle
K_{m,n}^*K_{0,m}^*f,g\rangle_{\mu_n}.
\]
Moreover, if $\mu_m$ is thought of as the expectation operator
$\mu_m\dvtx \ell^2(\mu_m)\ra\ell^2(\mu_n)$, $f\mapsto\mu_m(f)$,
then $(K_{m,n}^*-\mu_m)^*=K_{m,n}-\mu_n$.
Let
\[
\delta_x(z)=
\cases{
\mu_0(x)^{-1}, &\quad if $z=x$, \cr
0, &\quad otherwise.}
\]
Set $q=q(m)=2\prod_{i=1}^m(1+l(P_i))$ and $q'(m)$ to be the conjugate
exponent of
$q(m)$ so that $1/q(m)+1/q'(m)=1$. By duality, we have
\begin{eqnarray*}
&&
d_2(K_{0,n}(x,\cdot),\mu_n)\\
&&\qquad=\biggl\|\frac{K_{0,n}(x,\cdot)}{\mu_n(\cdot)}-1\biggr\|_{\ell^2(\mu
_n)} =
\biggl\|\frac{K_{0,n}^*(\cdot,x)}{\mu_0(x)}-1\biggr\|_{\ell^2(\mu_n)} \\
&&\qquad=\|(K_{0,n}^*-\mu_0)\delta_x\|_{\ell^2(\mu_n)}
=\|(K_{m,n}^*-\mu_m)K_{0,m}^*\delta_x\|_{\ell^2(\mu_n)}\\
&&\qquad\leq\|K_{m,n}^*-\mu_m\|_{\ell^2(\mu_m)\ra\ell^2(\mu_n)}\|
K_{0,m}^*\delta_x\|_{\ell^2(\mu_m)}\\
&&\qquad\leq\|\delta_x\|_{\ell^{q'(m)}(\mu_0)}\|K_{0,m}^*\|_{\ell^{q'(m)}(\mu
_0)\ra\ell^2(\mu_m)}
\|K_{m,n}^*-\mu_m\|_{\ell^2(\mu_m)\ra\ell^2(\mu_n)}\\
&&\qquad\leq\mu_0(x)^{-1/q(m)}\|K_{0,m}\|_{\ell^2(\mu_m)\ra\ell^{q(m)}(\mu_0)}
\|K_{m,n}-\mu_n\|_{\ell^2(\mu_n)\ra\ell^2(\mu_m)}.
\end{eqnarray*}
By assumption, we have that $q(m)\geq\log(\mu_0(x)^{-1})$, it now
follows from Corollary~\ref{cor-Sob} that
\[
d_2(K_{0,n}(x,\cdot),\mu_n)\leq e \prod_{i=m+1}^n\sigma_1(K_i,\mu_{i-1}).
\]
\upqed\end{pf}

%s4.2 ###
\subsection{Logarithmic Sobolev inequalities and $c$-stability}

\begin{theo}\label{thm-LS-c-stab}
Fix $c\in(1,\infty)$. Let $V$ be a finite set equipped with
a sequence of irreducible Markov kernels, $(K_i)_1^{\infty}$.
Assume that $(K_i)_1^{\infty}$ is $c$-stable with respect to a positive
probability measure $\mu_0$. For each $i$, set $\mu_0^i=\mu_0K_i$ and let
$\sigma_1(K_i,\mu_0)$ be the second largest singular value of the
operator $K_i\dvtx\ell^2(\mu_0^i)\ra\ell^2(\mu_0)$ and $l(K_i^*K_i)$ the
logarithmic Sobolev constant for the operator
$K_i^*K_i\dvtx\ell^2(\mu_0^i)\ra\ell^2(\mu_0^i)$. If
\[
\tilde{m}_x=\min\Biggl\{t\in\mathbb{N}\dvtx
\sum_{i=1}^t\log\bigl(1+c^{-2}l(K_i^*K_i)\bigr)\geq\log\log(\mu_0(x)^{-1/2})
\Biggr\},
\]
then for $n\geq\tilde{m}_x$ we have that
\[
d_2(K_{0,n}(x,\cdot),\mu_n)^2\leq
e^{2}\prod_{i=\tilde m_x+1}^n\bigl(1-c^{-2}\bigl(1-\sigma_1(K_i,\mu_0^i)^2\bigr)\bigr).
\]
\end{theo}
\begin{pf}
First, we note that $\mu_i/\mu_0^i\in[c^{-1},c]$. Let $P_i$ be the
Markov kernel
described in the statement of Theorem \ref{thm-LS}. By the same
arguments as
in Theorem~\ref{th-nash-stab}, we get that for all $x,y\in V$
\[
\mu_i(x)P_i(x,y)=\sum_z\mu_{i-1}(z)K_i(z,x)K_i(z,y)
\geq c^{-1}\mu_0^i(x)K_i^*K_i(x,y).
\]
A simple comparison argument similar to those used in the proof of
Theorem \ref{th-nash-stab} (see also \cite{DS-C,DS-L}) yields that
\[
l(P_i)\geq c^{-2}l(K_i^*K_i) \quad\mbox{and}\quad
1-\sigma(K_i,\mu_{i-1})^2\geq c^{-2}\bigl(1-\sigma(K_i,\mu_0^i)^2\bigr).
\]
The first inequality implies that $\tilde{m}_x\geq m_x$ where
$m_x$ is defined in the proof of Theorem \ref{thm-LS}. Using the results
of Theorem \ref{thm-LS} and the second inequality above gives the
desired result.
\end{pf}

The next result is when we have a $c$-stability assumption on a family
of kernels.
\begin{theo}
Let $c\in(1,\infty)$. Let $\mathcal Q$ be a family of irreducible
aperiodic Markov
kernels on a finite set $V$. Assume that $\mathcal Q$ is $c$-stable
with respect to some
positive probability measure $\mu_0$. Let $(K_i)_1^{\infty}$ be a
sequence of
Markov kernels with $K_i\in\mathcal Q$ for all $i$. Let $\pi_i$ be the
invariant measure of $K_i$. Let $\sigma_i(K_i)$ be the second largest singular
value for the operator $K_i\dvtx\ell^2(\pi)\ra\ell^2(\pi)$. Let $l(K_i^*K_i)$
be the logarithmic Sobolev constant for the operator $K_i^*K_i$ where
$K_i^*$ is the adjoint of $K_i\dvtx\ell^2(\pi)\ra\ell^2(\pi)$. If
\[
\tilde{m}_x=\min\Biggl\{t\in\mathbb{N}\dvtx
\sum_{i=1}^t\log\bigl(1+c^{-4}l(K_i^*K_i)\bigr)\geq\log\log(\mu_0(x)^{-1/2})
\Biggr\},
\]
then for $n\geq\tilde{m}_x$ we have that
\[
d_2(K_{0,n}(x,\cdot),\mu_n)^2\leq
e^{2}\prod_{i=m_x+1}^n\bigl(1-c^{-4}\bigl(1-\sigma_1(K_i)^2\bigr)\bigr).
\]
\end{theo}
\begin{pf}
Let $\mu_i=\mu_0K_{0,i}$. If $\mathcal Q$ is $c$-stable, then
$\mu_i/\pi_i\in[c^{-2},c^2]$. Similar arguments to those used in
Theorem \ref{thm-LS-c-stab} give the desired result.
\end{pf}

%s4.3 ###
\subsection{The relative sup norm}

To control the relative-sup merging time by this method,
we need an additional hypothesis. In the case of the $\ell^2$ distance,
we only required a control over the logarithmic Sobolev constant of
the kernel $P_i=K_i^*K_i\dvtx\ell^2(\mu_i)\ra\ell^2(\mu_i)$. In this
case, we will also need to control the logarithmic Sobolev constant
of $\check{P}_i=K_iK_i^*\dvtx\ell^2(\mu_{i-1})\ra\ell^2(\mu_{i-1})$ where
$K_i^*$ is the adjoint of the operator $K_i$ from $\ell^2(\mu_i)$ to
$\ell^2(\mu_{i-1})$.
\begin{theo}\label{thm-LS-sup}
Let $V$ be a finite set equipped with a sequence of Markov kernels
$(K_n)_0^{\infty}$
and an initial distribution $\mu_0$. Let $\mu_n=\mu_0K_{0,n}$ and
$P_i=K_i^*K_i\dvtx\ell^2(\mu_i)\ra\ell^2(\mu_i)$ and
$\check{P}_i=K_iK_i^*\dvtx\ell^2(\mu_{i-1})\ra\ell^2(\mu_{i-1})$
where $K_i^*$ is the adjoint of $K_i$ with respect to the measure $\mu_i$.
Let $l(P_i)$ and $l(\check{P}_i)$ be the logarithmic Sobolev constants of
$P_i$ and $\check{P}_i$, respectively. If $\mu_i^{\#}=\min_x\{\mu_i(x)\}
$ and
\begin{eqnarray*}
m_0^{\#}&=&\min\Biggl\{t\in\mathbb{N}\dvtx
\sum_{i=1}^t\log\bigl(1+l(P_i)\bigr)\geq\log\log({\mu_0^{\#}}{}^{-1/2}
)\Biggr\},\\
m_n^{\#}&=&\min\Biggl\{t\in\mathbb{N}\dvtx
\sum_{i=n-t}^n\log\bigl(1+l(\check{P}_i)\bigr)\geq\log\log({\mu_n^{\#
}}{}^{-1/2})\Biggr\},
\end{eqnarray*}
then for any $n\geq2m$,
\[
\max_{x,y}\biggl\{\biggl|\frac{K_{0,n}(x,y)}{\mu_n(y)}-1\biggr|\biggr\}
\leq e^2 \prod_{i=m+1}^{n-m}\sigma_1(K_i,\mu_{i-1}),
\]
where $m=\max\{m_0^{\#},m_n^{\#}\}$.
\end{theo}
\begin{rem}
This innocent looking theorem is not easy to apply. For instance, $m$
depends on $n$ and
without some control on this dependence the result is useless.
\end{rem}
\begin{pf*}{Proof of Theorem \ref{thm-LS-sup}}
Write
\[
\max_{x,y}\biggl\{\biggl|\frac{K_{0,n}(x,y)}{\mu_n(y)}-1\biggr|\biggr\}=
\|K_{0,n}-\mu_n\|_{\ell^1(\mu_n)\ra\ell^{\infty}(\mu_0)}
\]
and
\begin{eqnarray*}
&&\|K_{0,n}-\mu_n\|_{\ell^1(\mu_n)\ra\ell^\infty(\mu_0)}\\
&&\qquad\le
\|K_{n-m,n}-\mu_n\|_{\ell^1(\mu_n)\ra\ell^2(\mu_{n-m})}
\times
\|K_{m,n-m}-\mu_{n-m}\|_{\ell^2(\mu_{n-m})\ra\ell^2(\mu_m)}\\
&&\qquad\quad{}\times\|K_{0,m}-\mu_m\|_{\ell^2(\mu_m)\ra\ell^{\infty}(\mu_0)}.
\end{eqnarray*}
Note that
\[
\|K_{m,n-m}-\mu_{n-m}\|_{\ell^2(\mu_{n-m})\ra\ell^2(\mu_m)}
\leq\prod_{i=m+1}^{n-m}\sigma_1(K_{i,\mu_{i-1}})
\]
so we just need to bound the remaining terms in the right-hand side of the
inequality above. To bound $\|K_{n-m,n}-\mu_n\|_{\ell^1(\mu_n)\ra\ell
^2(\mu_{n-m})}$
set $q^*=q^*(m)=2\prod_{i=1}^m(1+l(\check{P}_{n-m+i}))$ and write
\begin{eqnarray*}
&&\|K_{n-m,n}-\mu_n\|_{\ell^1(\mu_n)\ra\ell^2(\mu_{n-m})}\\
&&\qquad=\|K_{n-m,n}^*-\mu_{n-m}\|_{\ell^2(\mu_{n-m})\ra\ell^{\infty}(\mu
_{n})}\\
&&\qquad =\|I(K_{n-m,n}^*-\mu_{n-m})\|_{\ell^2(\mu_{n-m})\ra\ell
^{\infty}(\mu_{n})}\\
&&\qquad \leq\|K_{n-m,n}^*-\mu_{n-m}\|_{\ell^2(\mu_{n-m})\ra\ell
^{q^*}(\mu_{n})}
\|I\|_{\ell^{q^*}(\mu_n)\ra\ell^{\infty}(\mu_n)}\\
&&\qquad \leq\|K_{n-m,n}^*\|_{\ell^2(\mu_{n-m})\ra\ell^{q^*}(\mu_{n})}
\|I\|_{\ell^{q^*}(\mu_n)\ra\ell^{\infty}(\mu_n)}.
\end{eqnarray*}
It follows from Corollary \ref{cor-Sob} that
\[
\|K_{n-m,n}-\mu_n\|_{\ell^1(\mu_n)\ra\ell^2(\mu_{n-m})}\leq
\|I\|_{\ell^{q^*}(\mu_n)\ra\ell^{\infty}(\mu_n)}\leq{\mu_n^{\#}}{}^{-1/q^*}.
\]
By assumption, we have that $q^*=q^*(m)\geq\log({\mu_n^{\#}}{}^{-1})$ so
we get
\[
\|K_{n-m,n}-\mu_n\|_{\ell^1(\mu_n)\ra\ell^2(\mu_{n-m})}\leq e.
\]
To bound $\|K_{0,m}-\mu_0\|_{\ell^1(\mu_m)\ra\ell^2(\mu_0)}$
set $q=q(m)=2\prod_{i=1}^m(1+l(P_i))$ and write
\[
\|K_{0,m}-\mu_0\|_{\ell^2(\mu_m)\ra\ell^{\infty}(\mu_0)}
\leq\|K_{0,m}-\mu_0\|_{\ell^2(\mu_{m})\ra\ell^{q}(\mu_0)}
\|I\|_{\ell^q(\mu_0)\ra\ell^{\infty}(\mu_0)}.
\]
It follows from Corollary \ref{cor-Sob} that
\[
\|K_{0,m}-\mu_0\|_{\ell^2(\mu_m)\ra\ell^{\infty}(\mu_0)}\leq
\|I\|_{\ell^{q}(\mu_0)\ra\ell^{\infty}(\mu_0)}\leq{\mu_0^{\#}}{}^{-1/q}.
\]
Since $q=q(m)\geq\log({\mu_0^{\#}}{}^{-1})$ we get
$\|K_{n-m,n}-\mu_n\|_{\ell^1(\mu_n)\ra\ell^2(\mu_{n-m})}\leq e$.
\end{pf*}
\begin{theo}
Fix $c\in(1,\infty)$.
Let $V$ be a finite set equipped with a sequence of Markov kernels
$(K_n)_1^{\infty}$.
Assume that $(K_n)_1^{\infty}$ is $c$-stable with respect to a positive
probability measure $\mu_0$. For each $i$, set $\mu_0^i=\mu_0K_i$
and $\mu_n^i=\mu_nK_i$. Let $\sigma_1(K_i,\mu_0^i)$ be the second largest
singular value of the operator $K_i\dvtx\ell^2(\mu_0^i)\ra\ell^2(\mu_0)$.
Let $l(K_i^*K_i)$ be the logarithmic Sobolev constant
of the operator $K_i^*K_i\dvtx\ell^2(\mu_0^i)\ra\ell^2(\mu_0^i)$ where
$K_i^*$ is the adjoint of $K_i\dvtx\ell^2(\mu_0^i)\ra\ell^2(\mu_0)$. Let
$l(K_iK_i^*)$ be the logarithmic Sobolev constant of the
operator $K_iK_i^*\dvtx\ell^2(\mu_n)\ra\ell^2(\mu_n)$ where
$K_i^*$ is the adjoint of $K_i\dvtx\ell^2(\mu_n^i)\ra\ell^2(\mu_n)$.
If $\mu_i^{\#}=\min_x\{\mu_i(x)\}$ and
\begin{eqnarray*}
\tilde{m}_0^{\#}&=&\min\Biggl\{t\in\mathbb{N}\dvtx
\sum_{i=1}^t\log\bigl(1+c^{-2}l(K_i^*K_i)\bigr)\geq\log\log({\mu_0^{\#
}}{}^{-1/2})\Biggr\},\\
\tilde{m}_n^{\#}&=&\min\Biggl\{t\in\mathbb{N}\dvtx
\sum_{i=n-t}^n\log\bigl(1+c^{-6}l(K_iK_i^*)\bigr)\geq\log\log({\mu_n^{\#
}}{}^{-1/2})\Biggr\},
\end{eqnarray*}
then for any $n\geq2\tilde m$
\[
\max_{x,y}\biggl\{\biggl|\frac{K_{0,n}(x,y)}{\mu_n(y)}-1\biggr|\biggr\}
\leq e^2\prod_{i=\tilde{m}}^{n-\tilde{m}}\bigl(1-c^{-2}\bigl(1-\sigma_1(K_i,\mu
_0^i)^2\bigr)\bigr)^{1/2},
\]
where $\tilde{m}=\max\{\tilde{m}_0^{\#},\tilde{m}_n^{\#}\}$.
\end{theo}
\begin{pf} Note that $\mu_i/\mu_0^i\in[c^{-1},c]$ and $\mu_i/\mu
_n^i\in[c^{-2},c^2]$.
Let $P_i$ and $\check{P}_i$ be the Markov kernels described in Theorem
\ref{thm-LS-sup} with kernels
%
%e4.4 ###
%e4.3 ###
%
\begin{eqnarray}
\label{eqn-P-LS-sup}
P_i(x,y)&=&\frac{1}{\mu_i(x)}\sum_z\mu_{i-1}(z)K_i(z,x)K_i(z,y),\\
\label{eqn-ch-P-LS-sup}
\check{P}_i(x,y)&=&\sum_z\frac{\mu_{i-1}(y)}{\mu
_i(z)}K_i(x,z)K_i(y,z).
\end{eqnarray}
Similar reasoning to that of Theorem \ref{thm-LS-c-stab} gives
\[
l(P_i)\geq c^{-2}l(K_i^*K_i)
\quad\mbox{and}\quad
1-\sigma(P_i)^2\geq c^{-2}\bigl(1-\sigma(K_i,\mu_0^i)^2\bigr),
\]
where $K_i^*$ above is the adjoint of $K_i\dvtx\ell^2(\mu_0^i)\ra\ell
^2(\mu_0)$.
This implies that $\tilde{m}_0^{\#}\geq m_0^{\#}$ where $m_0^{\#}$
is defined in Theorem \ref{thm-LS-sup}.

In the case of $\check{P}_i$, equation (\ref{eqn-ch-P-LS-sup}) gives
\[
\check{P}_i\geq c^{-4}\sum_z\frac{\mu_n(y)}{\mu_n^i(z)}K_i(x,z)K_i(y,z)
=c^{-4}K_iK_i^*(x,y),
\]
where $K_i^*$ is the adjoint of the operator $K_i\dvtx\ell^2(\mu
_n^i)\ra\ell
^2(\mu_n)$.
A simple comparison argument yields
\[
l(\check{P}_i)\geq c^{-6}l(K_iK_i^*)
\]
and so $\tilde{m}_n^{\#}\geq m_n^{\#}$ where $m_n^{\#}$
is defined in Theorem \ref{thm-LS-sup}. The desired result now
follows from Theorem \ref{thm-LS-sup}.
\end{pf}

The next theorem gives us similar results when we have $c$-stability
for a family of kernels.
\begin{theo}
Fix $c\in(1,\infty)$.
Let $\mathcal Q$ be a family of irreducible aperiodic Markov kernels on $V$.
Assume that $\mathcal Q$ is $c$-stable with respect to some positive probability
measure $\mu_0$. Let $(K_n)_1^{\infty}$ be a sequence with $K_i\in
\mathcal Q$
for all $i\geq1$. Let $\pi_i$ be the invariant measure of $K_i$ and
$\sigma_1(K_i)$
the second largest singular value of the operator $K_i$ acting on $\ell
^2(\pi_i)$.
Let $l(K_i^*K_i)$ and $l(K_iK_i^*)$ be the logarithmic Sobolev constants
of the operators $K_i^*K_i$ and $K_iK_i^*$ where $K^*_i$ is the adjoint of
$K_i\dvtx\ell^2(\pi_i)\ra\ell^2(\pi_i)$.
If $\mu_i^{\#}=\min_x\{\mu_i(x)\}$ and
\begin{eqnarray*}
\tilde{m}_0^{\#}&=&\min\Biggl\{t\in\mathbb{N}\dvtx
\sum_{i=1}^t\log\bigl(1+c^{-4}l(K_i^*K_i)\bigr)\geq\log\log({\mu_0^{\#
}}{}^{-1/2})\Biggr\},\\
\tilde{m}_n^{\#}&=&\min\Biggl\{t\in\mathbb{N}\dvtx
\sum_{i=n-t}^n\log\bigl(1+c^{-6}l(K_iK_i^*)\bigr)\geq\log\log({\mu_n^{\#
}}{}^{-1/2})\Biggr\},
\end{eqnarray*}
then for any $n\geq2\tilde m$
\[
\max_{x,y}\biggl\{\biggl|\frac{K_{0,n}(x,y)}{\mu_n(y)}-1\biggr|\biggr\}
\leq e^2
\prod_{i=\tilde{m}}^{n-\tilde{m}}\bigl(1-c^{-4}\bigl(1-\sigma_1(K_i)^2\bigr)\bigr)^{1/2},
\]
where $\tilde{m}=\max\{\tilde{m}_0^{\#},\tilde{m}_n^{\#}\}$.
\end{theo}
\begin{pf} First, note that $\mu_i/\pi_i\in[c^{-2},c^2]$. Equation
(\ref{eqn-P-LS-sup}) implies that
\[
l(P_i)\geq c^{-4}l(K_i^*K_i)
\quad\mbox{and}\quad
1-\sigma(K_i,\mu_i)^2\geq c^{-4}\bigl(1-\sigma(K_i)^2\bigr).
\]
To bound $l(\check{P}_i)$, we use (\ref{eqn-ch-P-LS-sup}) to get that
for all $x,y\in V$
\[
\check{P}_i(x,y)\geq c^{-4} K_iK_i^*(x,y).
\]
This implies that $l(\check{P}_i)\geq c^{-6}l(K_iK_i^*)$. It follows that
$\tilde{m}\geq m$ where $m$ is defined in Theorem \ref{thm-LS-sup}.
Applying Theorem \ref{thm-LS-sup} now gives us the desired result.
\end{pf}

%%%%%%%%%%%%%%%%%%%%%%%%%%%%%%%%%%%%%%%%%%%%%%%%%%%%%%%%%%%%%%%%%%%%%%%%%%%%%%%%%%%%%%%%%%%
%s4.4 ###
\subsection{An inhomogeneous walk on the hypercube}\label{subsec-hypercube}

Denote by $V=\{0,1\}^{2N}$ the $2N$-dimensional
hypercube, we say that $x,y\in V$ are neighbors, or $x\sim y$ if
\[
\sum_{i=1}^N|x_i-y_i|=1,
\]
where $x_i$ is the $i$th coordinate of $x\in V$. The simple
random walk
on $V$ is driven by the kernel
\[
K(x,y)=
\cases{
\dfrac{1}{2N}, &\quad if $x\sim y$, \vspace*{2pt}\cr
0, &\quad otherwise.}
\]
It is easy to check that $\mu$, the uniform measure on $V$, is
stationary for $K$.

Fix $\varepsilon\in(0,1)$ and consider the following perturbed version
of $K$.
\[
K_{\varepsilon}(x,y)=
\cases{
\dfrac{1}{2N}, &\quad if $x\sim y$ and $|x|\neq N$, \vspace*{2pt}\cr
\dfrac{1+\varepsilon}{2N}, &\quad if $x\sim y$ and $|x|=N, y=|N|+1$, \vspace*{2pt}\cr
\dfrac{1-\varepsilon}{2N}, &\quad if $x\sim y$ and $|x|=N, y=|N|-1$, \vspace*{2pt}\cr
0, &\quad otherwise.}
\]
For $\varepsilon\in(0,1)$, set
\[
\mathcal Q(\varepsilon)=\{K_{\delta}\dvtx \delta\in[-\varepsilon
,\varepsilon] \}.
\]
The example of time inhomogeneous Markov chains associated to
$Q(\varepsilon)$
above is related to the binomial example in \cite{SZ3}. See Remark \ref
{rem-bino} below.

We shall show that $\mathcal Q(\varepsilon)$ is $c$-stable. First,
consider the following
definition.
\begin{defin}
Let $\mathcal S_{2N}$ be the set of probability measures on $V= \{0$, $1\}
^{2N}$ that
satisfy the following three properties:
\begin{enumerate}[(2)]
\item[(1)] For all $x\in V$ with $|x|=N$ we have
$\nu(x)=\frac{1}{4^N}$.\vspace*{2pt}
\item[(2)] For all $i\in\{-N,\ldots,-1,1,\ldots,N\}$ there exists constants
$a_{\nu,i}$ such that $a_{\nu,i}=-a_{\nu,-i}$ and for any $x$ with
$|x|=N+i$ we have
\[
\nu(x)=\frac{1}{4^N}+a_{\nu,i}.
\]
\item[(3)] For all $i\in\{-N,\ldots,-1,1,\ldots,N\}$ we have $|a_{\nu
,i}|\leq\varepsilon/4^N$.
\end{enumerate}
\end{defin}
\begin{claim}\label{clm-hyper}
Let $\nu$ be in $\mathcal S_{2N}$ defined above, then for any $K\in
\mathcal Q(\varepsilon)$
we have that $\nu K\in\mathcal S_{2N}$.
\end{claim}
\begin{pf}
Let $\nu\in\mathcal S_{2N}$ and $Q\in\mathcal Q(\varepsilon)$, then
$Q=K_{\delta}$ for some
$\delta\in[-\varepsilon,\varepsilon]$. We will check each condition needed
for $\nu Q$
to be in $\mathcal S_{2N}$ separately.

(1)
For any $x$ with $|x|=N$ we have that $\nu Q(x)=\nu K(x)$.
The desired result now follows from the definition of $\mathcal S_{2N}$.

(2)
For $i$ such that $|i|\notin\{1,N\}$, consider an element $x$ such that
$|x|=N+i$.
Then
\begin{eqnarray*}
\nu Q(x)
&=&\mathop{\sum_{y\sim x}}_{|y|=|x|+1}\nu(y)Q(y,x)+
\mathop{\sum_{y\sim x}}_{|y|=|x|-1}\nu(y)Q(y,x)\\
&=&\biggl(\frac{1}{2N}\biggr)\biggl(
\mathop{\sum_{y\sim x}}_{|y|=|x|+1}\nu(y)+
\mathop{\sum_{y\sim x}}_{|y|=|x|-1}\nu(y)
\biggr)\\
&=&\biggl(\frac{1}{2N}\biggr)\biggl[
\biggl(\frac{1}{4^N}+a_{\nu,i+1}\biggr)|x|+\biggl(\frac{1}{4^N}+a_{\nu
,i-1}\biggr)(2N-|x|)
\biggr]\\
&=&\frac{1}{4^N}+\frac{1}{2N}\bigl(a_{\nu,i+1}|x|+a_{\nu
,i-1}(2N-|x|)\bigr).
\end{eqnarray*}
A similar computation as above yields that for an element $x$ with
$|x|=N-i$ we
have
\[
\nu Q(x)=\frac{1}{4^N}-\frac{1}{2N}\bigl(a_{\nu,i+1}|x|+a_{\nu
,i-1}(2N-|x|)\bigr).
\]

When $i=N$, and $x$ is such that $|x|=N+i=2N$, we have
\begin{eqnarray*}
\nu Q(x)
&=&\mathop{\sum_{y\sim x}}_{|y|=2N-1}\nu(y)Q(y,x)\\
&=&\biggl(\frac{1}{2N}\biggr)\biggl(\frac{1}{4^N}+a_{\nu,N-1}\biggr)(2N)\\
&=&\frac{1}{4^N}+a_{\nu,N-1}.
\end{eqnarray*}
When $i=-N$, and $x$ is such that $|x|=N-i=0$ we get
$\nu Q(x)=\frac{1}{4^N}-a_{\nu,N-1}$ as desired.

Finally, we check that cases for elements $x$ with $|x|=N\pm1$.
Consider an $x$ such that $|x|=N-1$, then
\begin{eqnarray*}
\nu Q(x)
&=&\mathop{\sum_{y\sim x}}_{|y|=N-2}\nu(y)Q(y,x)+
\mathop{\sum_{y\sim x}}_{|y|=N}\nu(y)Q(y,x)\\
&=&\biggl(\frac{1}{2N}\biggr)\biggl(\frac{1}{4^N}-a_{\nu,2}\biggr)(N-1)
+\biggl(\frac{1-\delta}{2N}\biggr)\biggl(\frac{1}{4^N}\biggr)(N+1)\\
&=&\frac{1}{4^N}-\frac{1}{2N}\biggl(a_{\nu,2}(N-1)+\frac{\delta
(N+1)}{4^N}\biggr).
\end{eqnarray*}
When $|x|=N+1$, then
\begin{eqnarray*}
\nu Q(x)
&=&\mathop{\sum_{y\sim x}}_{|y|=N}\nu(y)Q(y,x)+
\mathop{\sum_{y\sim x}}_{|y|=N+2}\nu(y)Q(y,x)\\
&=&\biggl(\frac{1+\delta}{2N}\biggr)\biggl(\frac{1}{4^N}\biggr)(N+1)
+\biggl(\frac{1}{2N}\biggr)\biggl(\frac{1}{4^N}+a_{\nu,2}\biggr)(N-1)\\
&=&\frac{1}{4^N}+\frac{1}{2N}\biggl(a_{\nu,2}(N-1)+\frac{\delta
(N+1)}{4^N}\biggr)
\end{eqnarray*}
as desired. We can now concluded that $a_{\nu Q,i}=-a_{\nu Q,-i}$.

(3)
From the calculations in part (2), we know that for $x$ with
$x=N+i$ and
$|i|\notin\{1,N\}$ and $|i|=N$ we have
\[
\nu Q(x)=\frac{1}{4^N}+\frac{1}{2N}\bigl(a_{\nu,i+1}|x|+a_{\nu
,i-1}(2N-|x|)\bigr)
\]
and
\[
\nu Q(x)=\frac{1}{4^N}+a_{\nu,N-1},
\]
respectively. It follows from the fact that for all $i$, $|a_{\nu
,i}|\leq\varepsilon/4^N$ that
for both cases above $|a_{\nu Q, i }|\leq\varepsilon/4^N$.
When $|i|=1$, we have that for $x$ with $|x|=N+i=N\pm1$
\begin{eqnarray*}
\nu Q(x)
&\leq&\frac{1}{4^N}+\frac{1}{2N}\biggl(a_{\nu,2}(N-1)+\frac{\varepsilon
(N+1)}{4^N}\biggr)\\
&\leq&\frac{1}{4^N}+\frac{1}{2N}\biggl(\frac{\varepsilon(N-1)}{4^N}+\frac
{\varepsilon(N+1)}{4^N}\biggr)\\
&=&\frac{1+\varepsilon}{4^N}.
\end{eqnarray*}
A similar calculation yields $\nu Q(x)\geq\frac{1-\varepsilon}{4^N}$. The
proof now follows from
the fact that $a_{\nu Q,i}=-a_{\nu Q,-i}$.
\end{pf}
\begin{claim}\label{cor-stab-hyper}
The set $\mathcal Q(\varepsilon)$ is $\frac{1+\varepsilon}{1-\varepsilon
}$-stable with
respect to any measure in $\mathcal S_{2N}$.
\end{claim}
\begin{pf}
Let $\mu_0\in\mathcal S_{2N}$. Let $(K_i)_1^{\infty}$
be any sequence of kernels such that $K_i\in\mathcal{Q}(\varepsilon)$
for all
$i\geq1$. Let $\mu_n=\mu_0K_{0,n}$, then by Claim \ref{clm-hyper}
we have that $\mu_n\in\mathcal S_{2N}$ and so for any $x\in V$
\[
\frac{1-\varepsilon}{1+\varepsilon}\leq\frac{\mu_n(x)}{\mu_0(x)}\leq
\frac
{1+\varepsilon}{1-\varepsilon}.
\]
\upqed\end{pf}

The kernels $K_{\delta}\in\mathcal Q(\varepsilon)$ drive periodic chains
that will
alternate between points with an even number of $1$'s and odd number of $1$'s.
So we will study following random walk driven by the kernel
\[
Q_\delta=\tfrac{1}{2}(I+K_{\delta}),
\]
where $I$ is the identity. Set
\[
\overline{\mathcal Q}(\varepsilon)=\{Q_{\delta}\dvtx\delta\in
[-\varepsilon
,\varepsilon]\}.
\]
\begin{claim}\label{clm-hyper-sing-LS}
Let $(K_i)_1^{\infty}$ be a sequence of Markov kernels such that
$K_i\in\overline{\mathcal Q}(\varepsilon)$
for all $i\geq1$. Let $\mu_0\in\mathcal S_{2N}$ be a positive measure,
and let
$\mu_n=\mu_0K_{0,n}$. Set $P_i=K_i^*K_i\dvtx\ell^2(\mu_i)\ra\ell^2(\mu_i)$
where $K_i^*$ is the
adjoint of $K_i\dvtx\ell^2(\mu_i)\ra\ell^2(\mu_{i-1})$. Let $\sigma
_1(K_i,\mu_i)$ and
be the second largest singular value of $K_i\dvtx\ell^2(\mu_i)\ra\ell
^2(\mu
_{i-1})$. Let
$l(P_i)$ be logarithmic Sobolev constant of $P_i$. Then
\[
\sigma_1(K_i,\mu_i)\leq1-C(\varepsilon)\frac{1}{2N}
\quad\mbox{and}\quad
l(P_i)\geq\frac{C(\varepsilon)}{4N},
\]
where $C(\varepsilon)=(1+\varepsilon)^{-2}(1-\varepsilon)^4$.
\end{claim}
\begin{pf}
Let $Q=2^{-1}(I+K_0)$ and $u$ be the uniform measure on $\{0,1\}^{2N}$.
Let $P_i(x,y)=K_i^*K_i\dvtx\ell^2(\mu_i)\ra\ell^2(\mu_i)$. Using the
$\frac{1+\varepsilon}{1-\varepsilon}$-stability of the sequence $(\mu
_n)_0^{\infty}$,
we get that
\begin{eqnarray*}
\mu_i(x)P_i(x,y)&=&\sum_z\mu_{i-1}(z)K_i(z,x)K_i(z,y)\\
&\geq&\frac{1-\varepsilon}{1+\varepsilon} \frac{u(x)}{u(x)}\sum_z
u(z)K_i(z,x)K_i(z,y)\\
&\geq&\frac{(1-\varepsilon)^3}{1+\varepsilon} \frac{u(x)}{u(x)}\sum_z u(z)
Q(z,x)Q(z,y)\\
&\geq&\frac{(1-\varepsilon)^3}{1+\varepsilon}u(x)Q^{(2)}(x,y).
\end{eqnarray*}
A simple comparison yields
\[
\mathcal E_{P_i,\mu_i}(f,f)\geq(1-\varepsilon)^3(1+\varepsilon)^{-1}
\mathcal E_{Q^{(2)},u}(f,f).
\]
Further comparison gives that
%
%e4.6 ###
%e4.5 ###
%
\begin{eqnarray}
\label{hyper-sigma}
1-\sigma_1(K_i,\mu_i)&\geq& C(\varepsilon)\bigl(1-\sigma_1(Q)\bigr),\\
\label{hyper-LS}
l(P_i)&\geq& C(\varepsilon)l\bigl(Q^{(2)}\bigr).
\end{eqnarray}
It is well known that for $K_0\dvtx\ell^2(u)\ra\ell^2(u)$ (the simple
random walk) we have
$2l(K_0)=1-\sigma_1(K_0)=1/N$. This implies that $\sigma_1(Q)=1-1/2N$.
The singular value inequality in Claim \ref{clm-hyper-sing-LS} now
follows from
(\ref{hyper-sigma}). For the rest of the proof, we note that
Lemma 2.5 of \cite{DS-N} tells us that
$\mathcal E_{Q^{(2)},u}(f,f)\geq\mathcal E_{Q,u}(f,f)$, and
so we get $l(Q^{(2)})\geq l(Q)$. The logarithmic Sobolev inequality now
follows from (\ref{hyper-LS}) and the fact that $l(Q)=1/4N$.
\end{pf}

By applying Theorem \ref{thm-LS} and Claim \ref{clm-hyper-sing-LS}, we
get the following
theorem.
\begin{theo}\label{th-hyper}
For any $\varepsilon\in(0,1)$ there exists a constant $D(\varepsilon)$ such
that the
total variation merging time of the sequence $(K_i)_1^{\infty}$ with
$K_i\in\overline{\mathcal Q}(\varepsilon)$ for all $i\in\{1,2,\ldots\}$ is
bounded by
\[
T_{\mathrm{TV}}(\eta)\leq D(\varepsilon)N(\log N+\log_+1/\eta).
\]
Moreover, we can chose $D(\varepsilon)$ such that
\[
\forall n\geq D(\varepsilon)N(\log{N}+\log_{+}1/\eta)\qquad
\max_{x,y,z\in V}\biggl\{\biggl|\frac{K_{0,n}(x,z)}{K_{0,n}(y,z)}-1
\biggr|\biggr\}\leq\eta.
\]
\end{theo}

We note that the relative-sup merging time bound is obtained with the
same arguments as those used at the end of the proof of Theorem \ref{thm-nash2}.
\begin{rem}\label{rem-bino}
The theorem above is closely related to the example in Section 5.2 of
\cite{SZ3} which studies a time inhomogeneous chain on $\{-N,\ldots, N\}
$ resulting from perturbations of a birth and death chain with binomial
stationary distribution. Both \cite{SZ3} and Theorem \ref{th-hyper}
give the correct upper bound on the merging time yet \cite{SZ3}
requires knowledge about the entire spectrum of the operators driving
the chain while the theorem above uses logarithmic Sobolev techniques.
\end{rem}

%s4.5 ###
\subsection{Modified logarithmic Sobolev inequalities and entropy}
Let $\nu$ and $\mu>0$ be two probability measures on $V$.
Define the relative entropy between $\mu$ and $\nu$ as
\[
\operatorname{Ent}_{\mu}(\nu)=\sum_{x\in V}\mu(x)\log\biggl(\frac{\mu(x)}{\nu
(x)}\biggr).
\]
It is well known that $\sqrt{2}\|\mu-\nu\|_{\mathrm{TV}}\leq\sqrt
{\operatorname{Ent}_{\nu}(\mu)}$.
Let $(K_n)_0^{\infty}$ be a sequence of Markov kernels on $V$,
$\mu_0$ be some initial distribution on $V$ and $\mu_n=\mu_0K_{0,n}$.
It follows by the triangle inequality that for any $x,y\in V$
\[
\|K_{0,n}(x,\cdot)-K_{0,n}(y,\cdot)\|_{\mathrm{TV}}
\leq\sqrt{2}\max_{x\in V}\sqrt{\operatorname{Ent}_{\mu
_n}(K_{0,n}(x,\cdot))}.
\]

Let $\alpha=\alpha(K,\nu)$ be the largest constant such that for any
probability measure $\mu$
\[
\operatorname{Ent}_{\nu K}(\mu K)\leq(1-\alpha)\operatorname{Ent}_{\nu
}(\mu).
\]
Let $\mu'=\mu K$ and $K^*\dvtx\ell^2(\mu)\ra\ell^2(\mu')$ be the
adjoint of $K\dvtx\ell^2(\mu')\ra\ell^2(\mu)$. Set
\[
\check{P}=KK^*\dvtx\ell^2(\mu)\ra\ell^2(\mu).
\]
In \cite{DLM}, the contraction constant $\alpha$ is related
to the so-called modified logarithmic Sobolev constant
\[
l'(\check{P})=\inf
\biggl\{\frac{\mathcal E_{\mu, \check{P}}(f^2,\log(f^2))}{\mathcal
L(f^2,\mu)}
\dvtx\mathcal L(f^2,\mu)\neq0, f\neq\mbox{constant}\biggr\}.
\]
\begin{pro}[(\cite{DLM}, Proposition 5.1)]\label{pro-DLM}
There exists a universal constant $0<\rho<1$ such that for any Markov kernel
$K$ and any probability measure $\mu$,
\[
\rho l'(\check{P})\leq\alpha(K,\mu)\leq l'(\check{P}),
\]
where $\check{P}=KK^*$ and $K^*$ is the adjoint of the operator
$K\dvtx\ell^2(\mu')\ra\ell^2(\mu)$, $\mu'=\mu K$.
\end{pro}
\begin{pro}\label{pro-rho-LS}
Referring to the proposition above,
\[
\rho\geq\log2 \biggl(\frac{1-\log2}{2}\biggr).
\]
\end{pro}
\begin{pf}
The proof of Proposition 5.1 in \cite{DLM} uses the fact that there
exists some
$0<\tilde{\rho}<1$ such that for all $x\in[-1,\infty)$
\[
0\leq\varphi(x)\leq\tilde{\rho}{}^{-1}\varphi(x/2),
\]
where
\[
\varphi(x)=(1+x)\log(1+x)-x.
\]
Let $f(x)=\varphi(x)-(2/(1-\log2))\varphi(x/2)$. We will show that for
all $x\in[-1,\infty)$ then $f(x)\leq0$. By differentiating $f$ we get
\[
f'(x)=\log(1+x)-\biggl(\frac{1}{1-\log{2}}\biggr)\log\biggl(\frac
{2+x}{2}\biggr)
\]
and
\[
f'''(x)=\frac{4\log2(1+x)+x^2\log2-3-2x}{(1+x)^2(2+x)^2(1-\log2)}.
\]
In particular, for $x\in[-1,0]$ we have $f'''(x)\leq0$.
This along with the fact that
\[
f'(-0.9)\leq0,\qquad f'(-0.1)>0 \quad\mbox{and}\quad f'(0)=0
\]
implies that there exists only one $z\in(-1,0)$ such that $f'(z)=0$.
It follows that $f$ is decreasing on $[-1,z]$ and $f$ is increasing
on $[z,0]$. Since $f(-1)=f(0)=0$, then for $x\in[-1,0]$ we have that
$f(x)\leq0$.

For $x\in[0,\infty)$, we note that
\[
f''(x)=\frac{1}{1+x}-\frac{1}{(1-\log2)(2+x)}\leq0,
\]
which implies that $f'(x)\leq f'(0)=0$. The fact that $f(x)\leq f(0)=0$ implies
$\tilde{\rho}=2/(1-\log2)$. The desired result follows from the fact
that the proof of
Proposition~5.1 in \cite{DLM} shows that
\[
\alpha(K,\mu)\geq\tilde{\rho} \log(2)l'(\check{P}).
\]
\upqed\end{pf}

The results in \cite{DLM} allow us to study merging via
logarithmic Sobolev constants.
\begin{pro}\label{pro-entropy}
Let $V$ be a finite state space equipped with a sequence of Markov kernels
$(K_n)_1^{\infty}$ and an initial distribution $\mu_0$.
Let $\mu_n=\mu_0K_{0,n}$ and
$\check{P}_i=K_iK_i^*\dvtx\ell^2(\mu_{i-1})\ra\ell^2(\mu_{i-1})$ where
$K_i^*$ is the
adjoint of
$K_i\dvtx\ell^2(\mu_i)\ra\ell^2(\mu_{i-1})$. Set $\mu_0^*=\min_x\mu_0(x)$
then for any
$x,y\in V$
\[
\|K_{0,n}(x,\cdot)-K_{0,n}(y,\cdot)\|_{\mathrm{TV}}\leq
\sqrt{2}\log\biggl(\frac{1}{\mu_0^*}\biggr)^{1/2}
\prod_{i=1}^n\bigl(1-\rho l'(\check{P}_i)\bigr)^{1/2},
\]
where $\rho$ is given in \textup{Propositions} \ref{pro-DLM} \textup
{and} \ref{pro-rho-LS}.
\end{pro}
\begin{pf}
We note that for any $x,y\in V$
\[
\|K_{0,n}(x,\cdot)-K_{0,n}(y,\cdot)\|_{\mathrm{TV}}
\leq\sqrt{2}\max_{x,y}\sqrt{\operatorname{Ent}_{\mu_n}(K_{0,n}(x,\cdot))}.
\]
Proposition 5.1 in \cite{DLM} gives that
\[
\operatorname{Ent}_{\mu_n}(K_{0,n}(x,\cdot))\leq
\operatorname{Ent}_{\mu_0}(\delta_x)
\prod_{i=1}^n\bigl(1-\rho l'(\check{P}_i)\bigr).
\]
The desired result now follows from the fact that
\[
\operatorname{Ent}_{\mu_0}(\delta_x)
=\log\biggl(\frac{1}{\mu_0(x)}\biggr)\leq\log\biggl(\frac{1}{\mu
_0^*}\biggr).
\]
\upqed\end{pf}

%s4.6 ###
\subsection{Biased shuffles}
In this section, we present two examples where the modified logarithmic Sobolev
inequality technique yields the correct merging time while the regular
logarithmic
Sobolev inequality technique does not.
Let $V_n=S_n$ be the symmetric group equipped with the uniform
probability measure $u$. Let $\tilde{Q}_i$ be the the kernel of transpose
$i$ with random, that is,
\[
\tilde{Q}_i(x,y)=
\cases{
1/n, &\quad if $x^{-1}y=(i,j)$ for $j\in[1,n]$, \cr
0, &\quad otherwise.}
\]
Let $Q_i=2^{-1}(I+\tilde{Q}_i)$ be the associated lazy chain.
It is known that the lazy chain has a mixing time of
of $2n\log n$. More precisely,
\[
t\ge2n(\log n+c) \quad\Rightarrow\quad
\max_{x,y}\biggl\{\frac{Q_i^t(x,y)}{u(y)}-1\biggr\}\leq2e^{-2c}\qquad
\forall x\in S_n.
\]
See, for example, \cite{SZ2}. The results of \cite{Goel} show that the
modified logarithmic
Sobolev constant for $Q_i$ is bounded by
\[
\frac{1}{n-1}\geq l'(Q_i)\geq\frac{1}{4(n-1)}.
\]

Set $\mathcal Q=\{Q_i,i=1,\ldots,n\}$. Since all $Q_i$ are reversible
with respect to the uniform distribution $u$, the set $\mathcal Q$ is
$1$-stable with respect to $u$. Using the methods of
\cite{SZ} (see also \cite{Ga,MPS}), one can prove
that for any sequence $(K_i)_1^\infty$ with $K_i\in\mathcal Q$ for all
$i\geq1$
we have
\[
t\ge2n(\log n+c) \quad\Rightarrow\quad
\max_{x,y}\biggl\{\frac{K_{0,n}(x,y)}{u(y)}-1\biggr\}\le2 e^{-2c}
\qquad\forall x\in S_n.
\]
The inequality above is due to the fact that the $Q_i$ are driven by probability
measures so the $\ell^2$ distance bounds the $\ell^{\infty}$
distance and the eigenvectors in Theorem~3.2 of
\cite{SZ3} drop out to give
%
%e4.7 ###
%
\begin{equation}\label{eqn-toprand-MLS}
d_2(K_{0,t}(x,\cdot),u)^2\leq\sum_{i=1}^{n!-1}\prod_{j=1}^t\sigma_i(K_j)^2.
\end{equation}
One can then group the singular values in the equality above since
the $Q_i$'s are all are images of each other under some inner
automorphism of $S_n$ which implies $\sigma_j(Q_i)=\sigma_j(Q_k)$ for all
$i,j,k$. For a more detailed discussion, see \cite{SZ}.

We now consider two variants of this example that cannot be treated using
the singular values techniques of \cite{Ga,SZ,SZ3} or the
logarithmic Sobolev inequality technique of Sections \ref
{subsec-Hyper}--\ref{subsec-hypercube}
but where the modified logarithmic Sobolev inequality does yield a
successful analysis.
This technique can be applied to the
two examples in this section because of the following three reasons:
\begin{enumerate}[(2)]
\item[(1)] any sequence $(K_i)_1^{\infty}$ of interest can be shown to
be $c$-stable with respect to some well chosen initial distribution;
\item[(2)]all the kernels $K_i$ driving the time inhomogeneous
process are directly comparable to the $Q_i$'s and,
\item[(3)] due to $(1)$ and the laziness of the $Q_i$'s we can
successfully estimate
the modified logarithmic Sobolev constants $l'(Q_iQ_i^*)=l'(Q_i^{(2)})$
to be of order $1/n$.
\end{enumerate}

%s4.6.1 ###
\subsubsection{Symmetric perturbations in $S_n$}
For the first variant, fix $\varepsilon\in(0,1)$ and
consider the set $\mathcal Q^{\#}(\varepsilon)$
of all Markov kernels $K$ on $S_n$ such that:
\begin{enumerate}[(a)]
\item[(a)] $K(x,y)=K(y,x)$ (symmetry) and
\item[(b)] $\forall x,y$ we have $(1-\varepsilon) Q_i(x,y)\le K(x,y)\le
(1+\varepsilon)Q_i(x,y)$
for some $i\in\{1,\ldots,n\}$.
\end{enumerate}
Hence, $\mathcal Q^{\#}(\varepsilon)$
is the set of all symmetric edge perturbations of kernels in $\mathcal Q$.
As we require symmetry, the uniform distribution is invariant for all
the kernels in $\mathcal Q^{\#}(\varepsilon)$. Now, what can be said of the
merging properties of sequences $(K_i)_1^\infty$ with $K_i\in
\mathcal Q^{\#}(\varepsilon)$? Unlike $\mathcal Q$, the kernels in
$\mathcal Q^{\#}(\varepsilon)$ are not invariant under left multiplication
in $S_n$.
So the eigenvectors of Theorem 3.2 in \cite{SZ3} do not drop out,
and we only get
\[
d_2(K_{0,t}(x,\cdot),u)^2\leq n!\prod_{i=1}^t\sigma_1(K_i,u)^2.
\]
Singular value comparison yields
$\sigma_1(K_i,u)\le1-(1-\varepsilon)/(2n)$ which gives
\[
t\ge(1-\varepsilon)^{-1}n(n\log n+2c) \quad\Rightarrow\quad
d_2(K_{0,t}(x,\cdot),u)\leq e^{-c} \qquad\forall x\in S_n.
\]
This indicates merging after order $n^2\log n$ steps
instead of the expected order $n\log n$ steps.
For any sequence $(K_i)_1^\infty$ with $K_i\in\mathcal Q^{\#
}(\varepsilon)$
for all $i\geq1$ set $\check{P_i}=K_iK_i^*$ where $K_i^*$ is the
adjoint of the operator $K_i\dvtx\ell^2(u)\ra\ell^2(u)$. A simple comparison
argument gives
\[
\check{P}_i(x,y)\geq(1-\varepsilon)^2Q^2_j(x,y)
\]
for some\vspace*{-1pt} $j\in[1,n]$. Further comparison yields
$l'(\check{P}_i)\geq(1-\varepsilon)^2l'(Q_j^2)$.
Lem\-ma~2.5 of \cite{DS-N} implies that $l'(Q_j^2)\geq l'(Q_j)$ so
$l'(\check{P}_i)$ is
of order\vspace*{1pt} at least $1/n$. Hence, there exists some constant
$C(\varepsilon
)$ independent
of $n$ such that
\[
\|K_{0,t}(x,\cdot)-K_{0,t}(y,\cdot)\|_{\mathrm{TV}}\leq\sqrt{2
\log n!}
\bigl(1-C(\varepsilon)/n\bigr)^{t/2}.
\]

In particular, for some constant $D(\varepsilon)$ we get
$T_{\mathrm{TV}}(\eta)\le D(\varepsilon)n(\log n +\log_+1/\eta)$.
To obtain a result for the relative-sup norm,
one can use the (nonmodified) logarithmic Sobolev technique as
the modified logarithmic Sobolev technique only gives bounds
in total variation. It is known that the logarithmic Sobolev
constant for top to random is of order $1/(n\log n)$,
see \cite{LY},
leading to results that are off by a factor of $\log n$.
This technique yields the best available result,
\[
t\ge C(\varepsilon)n\bigl((\log n)^2+c\bigr)\quad \Rightarrow\quad
\max_{x,y,z}\biggl\{\biggl|\frac{K_{0,t}(x,z)}{K_{0,t}(y,z)}-1
\biggr|\biggr\}
\leq e^{-c}.
\]

%s4.6.2 ###
\subsubsection{Sticky permutations}

We now consider a second variation on the transpose cyclic to random example.
Let $\rho\in S_n$, $\delta\in(0,1-Q_1(\rho,\rho))$ and
% be the element defined by $\rho(i)=n-i+1$.
%In a deck of cards, this element reverses the order of the deck by
%making a new pile from an old pile by always moving the top card.
%the choice of the distinguished vertex $\rho$ is arbitrary here.
consider the Markov kernel
\[
K(x,y)=
\cases{
Q_1(x,y), &\quad if $x\neq\rho$, \cr
Q_1(x,y)+\delta, &\quad if $x=y=\rho$, \cr
Q_1(x,y)-\delta/(n-1), &\quad if $x=\rho$ and $x^{-1}y=(1,j)$ for $j\in[2,n]$.}
\]
In words, $K$ is obtained from $Q_1$ by adding extra holding
probability at $\rho$, making $\rho$ ``sticky.''
Next, if $\sigma$ is the cycle $(1,\ldots,n)$, let
\[
K_i(x,y)= K(\sigma^{i-1} x \sigma^{-i+1},\sigma^{i-1} y
\sigma^{-i+1}).
\]
In words, $K_i$ is $Q_i$ with some added holding at
$\rho_i=\sigma^{-i+1}\rho\sigma^{i-1}$.

We would like to consider the merging properties of the
sequence $(K_i)_1^\infty$. Unlike the previous example,
the uniform probability is not invariant under $K_i$.
However, this type of construction is considered in
\cite{SZ-wave}.

Let
\[
\varepsilon=\frac{\delta}{\sum_{z\neq\rho}Q_1(\rho,z)}
\]
so that $K(x,y)\geq(1-\varepsilon)Q_1(x,y)$.
It is proved that $(K_i)_1^\infty$ is $(1-\varepsilon)^{-1}$-stable
with respect
to the probability measure $\mu_0=\tilde{\pi}$, where
$\tilde{\pi}$ is the invariant probability measure of
the Markov kernel $\tilde{K}(x,y)=K(x,\sigma^{-1}y\sigma)$.
From the analysis in \cite{SZ-wave}, Section 5, one can see that
\[
(1-\varepsilon)u\leq\tilde{\pi}\leq(1-\varepsilon)^{-1} u.
\]
Applying the singular value techniques used in Section 5 of \cite{SZ-wave}
would give us an upper bound on the relative sup merging time of
order $n^2\log n$.

Set $\check{P}_i=K_iK_i^*\dvtx\ell^2(\mu_{i-1})\ra\ell^2(\mu_{i-1})$
where $K_i^*$ is the adjoint of the operator $K_i\dvtx\ell^2(\mu_i)\ra
\ell
^2(\mu_{i-1})$.
Since $K_i(x,y)\geq(1-\varepsilon)Q_i(x,y)$, for $x\neq y$ we can write
\[
\check{P_i}(x,y)=\sum_zK_i(x,z)K_i(y,z)\mu_{i-1}(y)\mu_i(z)^{-1}\geq
(1-\varepsilon)^4 Q_i^2(x,y).
\]
It follows by comparison that $l'(\check{P}_i)\geq(1-\varepsilon)^5
l'(Q_i^{2})$.
We can successfully estimate $l'(Q_i^{2})$ due to
Lemma 2.5 of \cite{DS-N} which implies $l'(Q_i^{2})\geq l'(Q_i)$. So we have
that $l'(\check{P}_i)$ is at least $(1-\varepsilon)^5/(4(n-1))$.
Proposition \ref{pro-entropy} gives us that
\[
\|K_{0,t}(x,\cdot)-K_{0,t}(y,\cdot)\|_{\mathrm{TV}}
\leq\sqrt{2} \log\biggl(\frac{n!}{1-\varepsilon}\biggr)^{1/2}
\biggl(1-\frac{\rho(1-\varepsilon)^5}{4(n-1)}\biggr)^{t/2},
\]
where $\rho$ is as in Proposition \ref{pro-rho-LS}.
So for some constant $D=D(\varepsilon)$, we get
\[
T_{\mathrm{TV}}(\eta)\le Dn\bigl(\log n+\log_+ (1/\eta)\bigr).
\]

% imsref loaded by lrinkeviciute, 2010-10-07 08:30:57
%

%
\printaddresses

\end{document}